\documentclass[11pt]{article}
\usepackage{amsfonts}
\usepackage{amssymb}
\usepackage{graphicx}
\usepackage{amsmath}
\usepackage{amsthm}
\usepackage{color}
\usepackage{multirow}
\usepackage{url}
\usepackage{float}
\usepackage{booktabs} 
\usepackage{makecell}
\usepackage{listings}
\usepackage{subcaption}
\usepackage{algorithm}
\usepackage{enumerate}
\usepackage{tikz}

\usepackage{algpseudocode}
\lstnewenvironment{Rcode}{\lstset{language=R}}{}

\newtheorem{theorem}{Theorem}[section]

\newtheorem{corollary}[theorem]{Corollary}

\newtheorem{lemma}[theorem]{Lemma}

\newtheorem{proposition}[theorem]{Proposition}
\newtheorem{remark}[theorem]{Remark}

\definecolor{battleshipgrey}{rgb}{0.52, 0.52, 0.51}
\definecolor{navyblue}{rgb}{0.0, 0.0, 0.5}
\definecolor{arsenic}{rgb}{0.23, 0.27, 0.29}
\definecolor{oldmauve}{rgb}{0.4, 0.19, 0.28}
\usepackage[colorlinks=TRUE]{hyperref}
\hypersetup{
	colorlinks=true,
	linkcolor=navyblue,
	filecolor=blue,      
	urlcolor=oldmauve,
	citecolor=navyblue,
	pdftitle={Overleaf Example},
	pdfpagemode=FullScreen,
}
\usepackage{natbib}
\setcitestyle{authoryear}

\begin{document}

	\title{ \bf Automatic and location-adaptive estimation in functional single-index regression}
	\author{Silvia Novo{$^a$}\footnote{Corresponding author email address: \href{s.novo@udc.es}{s.novo@udc.es}} \hspace{2pt} Germ\'{a}n Aneiros{$^b$} \hspace{2pt} Philippe Vieu{$^c$} \\		
		{\normalsize $^a$ Department of Mathematics, MODES, CITIC, Universidade da Coruña, A Coruña, Spain}\\
		{\normalsize $^b$ Department of Mathematics, MODES, CITIC, ITMATI, Universidade da Coruña, A Coruña, Spain}\\
		{\normalsize $^c$ Institut de Math\'{e}matiques, Universit\'e  Paul Sabatier, Toulouse, France}
	}
	
	\date{}
	\maketitle
	\begin{abstract} This paper develops a new automatic and location-adaptive procedure for estimating regression in a Functional Single-Index Model (FSIM). This procedure is based on $k$-Nearest Neighbours ($k$NN) ideas. The asymptotic study includes results for automatically data-driven selected number of neighbours, making the procedure directly usable in practice. The local feature of the $k$NN approach insures higher predictive power compared with usual kernel estimates, as illustrated in some finite sample analysis. As by-product we state as preliminary tools some new uniform asymptotic results for kernel estimates in the FSIM model.
	\end{abstract}
	
	
\section{Introduction}
\label{intr}
In regression analysis, one of the key questions is to construct methods able to balance the trade-off between too much flexibility of the model and easiness of implementation of the statistical procedure. Because of sparseness of data in high-dimensional spaces this need is more and more crucial when the dimensionality of the problem increases. In multivariate regression analysis, this has been the starting point for many advances around semiparametric modelling (see \citealt*{harmsw04}). In the functional data framework, regression problems involve infinite-dimensional variables and the dimensionality problem becomes even more important  (\Citealt{gee11}; \Citealt{vie17}).

Functional Data Analysis (FDA) has been really popularized through the works by Ramsay and Silverman (see \Citealt*{rams05}), and there is an extensive literature on functional regression modelling (see \Citealt{gres17} for a general presentation). This literature is mainly concentrated either around nonparametric models which were popularized by \cite{fer2006} (see \citealt*{gee11}, \citealt*{linv18} for recent surveys) or around linear models (see Chapter 11 of \citealt*{hsie15}), but the semiparametric framework is still a very underdeveloped field in FDA (see however \citealt*{goiv14} for an early discussion).

The recent advances in FDA have highlighted the necessity of developing models aiming to reduce dimensionality effects (see \citealt*{cue14}, \citealt*{goiv17}, \citealt*{vie17}, \citealt*{anecfv18} for recent surveys), and semiparametric ideas are natural candidates for that purpose. In that way, \citet*{ferpv03} and \citet*{ait} studied the Functional Single-Index Model (FSIM) 
\begin{equation}
	\label{modelo}
	Y=r\left(\left<\theta_0,X\right>\right)+\varepsilon,
\end{equation}
where $Y$ denotes a scalar response verifying some finite moment conditions (for details, see Section \ref{sec-assump-known}), $X$ is a functional explanatory random variable valued in a separable Hilbert space $\mathcal{H}$ with inner product $\left\langle \cdot, \cdot \right\rangle$, $\varepsilon$ is a random error verifying $\mathbb{E}\left(\varepsilon|X\right)=0$, $\theta_0 \in \mathcal{H}$ is the functional index and $r(\cdot)$ is the unknown link function. \cite{ferpv03} focused on the case of known $\theta_0$ and obtained the pointwise rate of convergence of a kernel estimator of $r\left(\left<\theta_0,x\right>\right)$, where $x \in \mathcal{H}$. The case of unknown $\theta_0$ was dealt in \cite{ait}, where both the consistency and optimality of a cross-validation based estimator of $\theta_0$ were proven. In addition, the FSIM (\ref{modelo}) was extended in different directions in \citet*{bouafv10}, \citet*{chehm11}, \citet*{ferr2013}, \citet*{ma16} and \citet*{wang2016}, among others (see Section \ref{tracks} for a slight discussion on such extensions and their differences with this paper).

This paper presents a wide study of the functional semiparametric model FSIM (\ref{modelo}). Section \ref{FSIM}   
develops a new automatic and location-adaptive procedure for estimating regression in the FSIM based on $k$-Nearest Neighbors ($k$NN) ideas (note that, in nonparametric statistics, an estimator is said to be ``location-adaptive'' when the smoothing parameter depends on the value in which one wishes to estimate, $x$; in the particular case of nonparametric regression estimation by means of the $k$NN estimator, the corresponding smoothing parameter is a bandwidth depending on the fixed value $k$ as well as on $x$; for details, see (\ref{est-k})). Section \ref{asymptotics} states general asymptotic results for the $k$NN procedure, with main interest of being uniform over all the parameters of the model. As discussed in Section \ref{data-driven}, this uniformity feature allows to derive directly results for random data-driven choices of these parameters making our procedure directly applicable in practice. Even if our main goal is to study $k$NN procedures, we derive also along Sections \ref{asymptotics} and \ref{data-driven} similar results for the standard kernel approach. The main feature of the obtained rates of convergence is that they are similar to those obtained in one-dimensional problems, giving evidence of the dimensionality reduction property of the method. Suggestions to address some practical issues related to the proposed methodology are shown in Section \ref{practical}. Such suggestions are supported in Section \ref{simulation} by means of a simulation study which, in addition, compares the performance of both the $k$NN-based and kernel-based procedures. Section \ref{application} illustrates, through some benchmark real curves dataset, how the  $k$NN approach outperforms standard procedures. It also shows how the semiparametric feature of the FSIM has not only nice predictive performance but it also provides easily interpretable and representable outputs. Finally, Section \ref{tracks} lists some tracks for future. The proofs of the main results are presented in the Appendix.

\section{The functional single-index model}\label{FSIM}
\subsection{Motivation} \label{motivation}
This paper focuses on the FSIM (\ref{modelo}), which 
can be seen as an extension of the standard well-known functional linear model (see \citealt*{hsie15} for discussion)
\begin{equation}
	\label{FLM}
	Y=\int_\mathcal{I}\theta_0(t)X(t)dt+\varepsilon ,
\end{equation}
as well as a special case of the functional nonparametric model (see \citealt*{fer2006})
\begin{equation}
	\label{FNPM}
	Y=m(X)+\varepsilon.
\end{equation}
\noindent In fact, the FSIM (\ref{modelo}) is an appealing trade-off between these two approaches. On the one hand, it is much more flexible, and hence more reliable in practice, than the parametric model (\ref{FLM}). On the other hand, it presents much less sensitivity to dimensionality effects than the nonparametric model (\ref{FNPM}) since it involves the estimation of the one-dimensional function $r$ rather than of the nonlinear infinite-dimensional operator $m$. These facts allow to say that the FSIM (\ref{modelo}) is a nice competitor for models (\ref{FLM}) and (\ref{FNPM}).

The model (\ref{modelo}) has been introduced in \cite{ferpv03} and conditions insuring its identifiability have been stated (see also \citealt*{ferpv11}). Here we assume that $r(\cdot)$ is differentiable and the following normalization constraint insures the uniqueness of the pair $(r,\theta_0)$:
\begin{equation}\label{identification}
	\left\langle\theta_0,e_1\right\rangle=1,
\end{equation}
\noindent  where $e_1$ is the first element of an orthonormal basis of $\mathcal{H}$ (for details, see \citealt{ferpv03}).

\subsection{The statistics}\label{model}
Let $Z_i = (X_i, Y_i)$, $i=1,...,n,$ be a sample of $n$ independent pairs identically distributed as $Z = (X,Y)$, which verifies the FSIM (\ref{modelo}); that is,
\begin{equation*}
	Y_i=r\left(\left<\theta_0,X_i\right>\right)+\varepsilon _i \ (i=1,\ldots,n).
\end{equation*}

\noindent For any $\theta \in \mathcal{H}$, we first consider the operator $$r_{\theta}(\cdot):\mathcal{H}\longrightarrow\mathbb{R}$$ defined as
$r_{\theta}(x)=r\left(\left<\theta,x\right>\right)$, $\forall x \in \mathcal{H}$, and we denote
$$d_{\theta}\left(x_1,x_2\right)=\left|\left<\theta,x_1-x_2\right>\right|, {\mbox{ for }} x_1,x_2\in\mathcal{H}.$$
\noindent For each direction $\theta$, we construct the $k$NN statistic as
\begin{equation}
	\hat{r}^\ast_{k,\theta}(x)=\frac{\sum_{i=1}^nY_iK\left(H_{k,x,\theta}^{-1}d_{\theta}(X_i,x)\right)}{\sum_{i=1}^nK\left(H_{k,x,\theta}^{-1}d_{\theta}(X_i,x)\right)}, \ \forall x\in\mathcal{H},
	\label{est-k}
\end{equation}
\noindent where $k\in \mathbb{Z}^+$ is a smoothing factor ($k=k_n$ depends on $n$) and $K$ is a kernel. In addition, we have denoted
\[
H_{k,x,\theta}=\min\left\{h\in \mathbb{R}^+ \mbox{\text{ such that }} \sum_{i=1}^n1_{B_\theta(x,h)}(X_i)=k\right\}
\]
with
\begin{equation*}
	B_{\theta}(x,h)=\left\{z\in \mathcal{H}:d_{\theta}\left(x,z\right) \leq h\right\}.
\end{equation*}
\noindent The $k$NN statistic $\hat{r}^\ast_{k,\theta}$ can be seen as an extension of the usual kernel statistic
\begin{equation}
	\hat{r}_{h,\theta}(x)=\frac{\sum_{i=1}^nY_iK\left(h^{-1}d_{\theta}(X_i,x)\right)}{\sum_{i=1}^nK\left(h^{-1}d_{\theta}(X_i,x)\right)}, \ \forall x\in\mathcal{H},
	\label{est}
\end{equation}
in which $h\in \mathbb{R}^+$ is the bandwidth ($h=h_n$ depends on $n$).

The $k$NN statistic is more appealing than the kernel one for two reasons. On the one hand, it involves a local smoothing factor $H_{k,x,\theta}$ making possible to capture local features of the data (while the smoothing factor $h$ of the kernel statistic does not depend on $x$). On the other hand, this local smoothing factor depends only on a discrete parameter $k$ taking values in the finite set $\{1,2,\ldots,n\}$. This fact makes much easier to select  $k$ in practice than the bandwidth $h$ appearing in kernel methods (which takes values in a continuous interval). In counterpart, the price to pay for so much flexibility of the procedure is that the theoretical properties are much more difficult to analyse (because of the randomness of the smoothing factor  $H_{k,x,\theta}$). More precisely, neither of the two terms in the ratio (\ref{est-k}) can be written as a sum of independent and identically distributed variables (as could be those appearing in (\ref{est})), and hence their analysis will require much more sophisticated tools than standard limit theorems for i.i.d. sequences. These features of $k$NN estimates have been widely highlighted in one-dimensional problems (see \citealt*{gyorfi} for a general discussion), but really few advances have been developed in the functional regression setting. The existing literature on kNN functionnal regression is mainly concerning nonparametric modelling (see \citealt*{biacg10}, \citealt*{kudv14}, \citealt*{muled}, \citealt*{kara} and \citealt*{linmv19} for the most recent advances, and see \citealt{linv18} for an exhaustive survey) or partial linear modelling (see \citealt*{linav17}), but at our knowledge this paper is stating the first advances in functional semiparametric regression.

In the next Sections \ref{asymptotics} and \ref{data-driven} we provide a complete study of the $k$NN procedure in the semiparametric model (\ref{modelo}). The main idea is to state asymptotic results in a uniform sense over all the parameters of the method (that is, over the direction $\theta$ and over the smoothing factor $k$). This will be done by following the uniform in bandwidth ideas widely developed in un-functional setting (see eg \citealt*{dony}) and recently adapted to functional setting (see \citealt*{kara2017a}), but including suitable technical changes to adapt such ideas both to $k$NN procedures and to the infinite-dimensional parameter $\theta$. Note that, even if our main goal is the study of the $k$NN procedure, we also derive as by-product a full asymptotic analysis of the standard kernel statistic (\ref{est}), extending in this case earlier results by \cite{ferpv03}, \cite{ait}, \cite{chehm11}, \cite{ferr2013}. 

\section{Asymptotic results} 
\label{asymptotics}
\subsection{Presentation and general notation}\label{notations}
Section \ref{as-1} starts by presenting the uniform in bandwidth (UIB) and uniform in the number of neighbours (UINN) consistency of the statistics $\hat{r}_{h,\theta}(x)$ (\ref{est}) and $\hat{r}^\ast_{k,\theta}(x)$ (\ref{est-k}), respectively, when $\theta$ is fixed. Then, Section \ref{as-2} extends these asymptotics by providing also uniform consistency over the functional parameter $\theta$.

Let us first introduce some notation. Throughout this paper, $x$ denotes a fixed element in $\mathcal{H}$ while $\theta$ is some direction in $\mathcal{H}$. Furthermore, we will use the notation:
\[
\phi_{x,\theta}(h)=\mathbb{P}\left(d_{\theta}(X, x) \leq h\right),
\]
\begin{equation}
	\mathcal{K}_{\theta} = \left\{\cdot\longrightarrow K\left(h^{-1}d_{\theta}(x,\cdot)\right), \ h > 0\right\} \label{h6}\end{equation}
and
\begin{equation}\mathcal{K}_{\Theta_n} = \cup_{\theta \in \Theta_n}\mathcal{K}_{\theta}=\left\{\cdot\longrightarrow K\left(h^{-1}d_{\theta}(x,\cdot)\right), \ h > 0,\ \theta\in\Theta_n\right\},\label{h666}\end{equation}
where $\Theta_n \subset \mathcal{H}$ is the set of directions of interest (note that both $\mathcal{K}_{\theta}$ and $\mathcal{K}_{\Theta_n}$ are classes of functions that should satisfy conditions (\ref{h7}) and (\ref{h777}), respectively; furthermore, condition (\ref{h11}) allows the size of $\Theta_n$ to grow up to infinite as $n$ does). In addition, let $\mathcal{Q}$ be a probability measure on the space $(\mathcal{H},\mathcal{A})$. Then, $||\cdot||_{\mathcal{Q},2}$ means the norm $L_2(\mathcal{Q})$ defined on certain space $S=\{f:\mathcal{H}\longrightarrow \mathbb{R}\}$, while $d_{\mathcal{Q},2}(\cdot,\cdot)$ is the metric associated to the norm $L_{2}(\mathcal{Q})$; that is, for $f,g \in S$, 

\begin{equation*}
	||f||_{\mathcal{Q},2}=\left(\int_{\mathcal{H}}\left|f(t)\right|^2 d\mathcal{Q}(t)\right)^\frac{1}{2}
\end{equation*}
and
\begin{equation*}
	d_{\mathcal{Q},2}(f,g)=||f-g||_{\mathcal{Q},2}=\left(\int_{\mathcal{H}}\left|f(t)-g(t)\right|^2 d\mathcal{Q}(t)\right)^\frac{1}{2}.
\end{equation*}
Finally, given a metric space $(\mathcal{K},d)$, $\mathcal{N}\left(\epsilon,\mathcal{K},d\right)$ denotes the minimal number of open balls (in the topological space given by $d$) with radius $\epsilon$ which are needed to cover $\mathcal{K}$.

\subsection{The case of $\theta_0$ known}\label{as-1}
Let us assume that the true direction, $\theta_0$, in the FSIM (\ref{modelo}) is known. In order to state the UIB and the UINN almost-complete
convergence of the  estimators $\hat{r}_{h,\theta_0}(x)$  and $\hat{r}^\ast_{k,\theta_0}(x)$ (\ref{est-k}), some of the following assumptions will be used. 

\subsubsection{Assumptions for UIB and UINN consistency}
\label{sec-assump-known}
\begin{description}		
	\item[\textit{About the small-ball probability.}] Let us assume that:
	\begin{itemize}
		\item For all $h>0$,
		\begin{equation}
			\phi_{x,\theta_0}(h)>0.
			\label{h1}
		\end{equation} 
		
		\item There exist a constant $0<C_1$ and sequences $\{a_n\},\{b_n\} \subset \mathbb{R}^+$ ($a_n \leq b_n$) such that, for $h \in [a_n,b_n]$ with $n$ large enough,
		\begin{equation}
			C_1 \leq\frac{\phi_{x,\theta_0}(h/2)}{\phi_{x,\theta_0}(h)}.
			\label{h2}
		\end{equation}
		
		\item The sequences $\{a_n\}$ and $\{b_n\}$ verify:
		\begin{equation}
			a_n \rightarrow 0, \ b_n \rightarrow 0 \mbox{ and } \frac{\log n}{n\min\left\{a_n,\phi_{x,\theta_0}(a_n)\right\}} \rightarrow 0.\label{h8}
		\end{equation}

		\item There exist sequences $\{\rho_n\} \subset (0,1)$, $\{k_{1,n}\}\subset \mathbb{Z}^+$ and $\{k_{2,n}\} \subset \mathbb{Z}^+$ ($k_{1,n}\leq k_{2,n}\leq n$) such that 
		\begin{equation}	
			\phi_{x,\theta_0}^{-1}\left(\frac{ k_{2,n}}{\rho_n n}\right) \rightarrow 0 ,
			\label{h33bis0}
		\end{equation}
		\begin{equation}
			\min \left\{\frac{1-\rho_n}{4}\frac{k_{1,n}}{\ln n},\frac{(1-\rho_n)^2}{4\rho_n}\frac{k_{1,n}}{\ln n}\right\} > 2 \label{hnew0}
		\end{equation}
		and
		\begin{equation}
			\frac{\log n}{n \min\left\{\phi_{x,\theta_0}^{-1}(\rho_n k_{1,n}/n), \rho_n k_{1,n}/n\right\}}\rightarrow 0 
			\label{h9}
		\end{equation}
		
	\end{itemize}

	\item[\textit{About the model.}] We assume that:
	\begin{itemize}
		\item There exist constants $\beta_0>0$ and $C_3>0$, such that:
		\begin{equation}
			\forall x_1, x_2\in N_{x,\theta_0}, \quad \left|r_{\theta_0}(x_1)-r_{\theta_0}(x_2)\right|\leq C_3d_{\theta_0}\left(x_1, x_2\right)^{\beta_0},
			\label{h3}
		\end{equation}
		where $N_{x,\theta_0}$ denotes a fixed neighbourhood of $x$ in the topological space induced by the semi-metric $d_{\theta_0}(\cdot,\cdot)$.
		
		\item 
		There exist constants $m\geq2$ and $C_4 > 0$, such that:
		\begin{equation}\mathbb{E}\left(|Y|^m|X\right)<C_4<\infty,\ a.s\label{h}.\end{equation}
	\end{itemize}
	
	\item[\textit{About the kernel.}] We assume that:
	\begin{itemize}
		\item  There exist constants $0< C_5 \leq C_6<\infty$, such that:
		\begin{equation}
			0<C_5 1_{(0,1/2)}(\cdot) \leq K(\cdot) \leq C_6 1_{(0,1/2)}(\cdot),\label{h4}
		\end{equation}
		where $1_{(0,1/2)}$ denotes the indicator function of the set $(0,1/2)$.
		\item The class of functions $\mathcal{K}_{\theta_0}$ (see (\ref{h6})) is a pointwise measurable class such that
		\begin{equation}
			\sup_{\mathcal{Q}}\int_{0}^1\sqrt{1+\log\mathcal{N}\left(\epsilon||F_{\theta_0}||_{\mathcal{Q},2},\mathcal{K}_{\theta_0},d_{\mathcal{Q},2}\right)}d\epsilon<\infty,
			\label{h7}
		\end{equation}
		where $F_{\theta_0}$ is the minimal envelope function of the set $\mathcal{K}_{\theta_0}$ and the supremum is taken  over all probability measures $\mathcal{Q}$ on the measurable  space  $(\mathcal{H},\mathcal{A})$ with $||F_{\theta_0}||_{\mathcal{Q},2}^2<\infty$.
	\end{itemize}

\end{description}

Assumptions (\ref{h1}), (\ref{h2}) and (\ref{h3})-(\ref{h4}) are standard ones in the setting of functional nonparametric regression models (see \citealt{fer2006}), while assumptions (\ref{h8}) and (\ref{h7}) are usual to obtain UIB consistency in such setting (see \citealt{kara2017a}). In fact, assumptions (\ref{h8}) and (\ref{h7}) adapt the ones used in that paper to the case where the nonparametric regression function is $r_{\theta_0}(\cdot)$ and the semi-metric to use in the kernel estimator is $d_{\theta_0}(\cdot,\cdot)$. Focusing now on the UINN consistency, assumptions (\ref{h33bis0})-(\ref{h9}) adapt (in the same way as in the previous case of UIB consistency) and correct those used in \cite{kara}. More specifically, in \cite{kara} it was forgotten to include the parameter $\alpha$ in their expression (17). As a consequence, Assumption (H4) in \cite{kara} should be modified in the way of our assumptions (\ref{h33bis0})-(\ref{h9}), where the notation $\rho_n$ was considered instead of $\alpha$; in addition, $\alpha$ should be introduced in the rates of convergence corresponding to their Theorem 3.1 in the same way as $\rho_n$ in our Theorem \ref{th1}(b). The justification of those changes in both the assumptions and the rates of convergence in \cite{kara} can be seen in the proof of our Theorem \ref{teorema theta estimado}(b). Finally, in the particular case of assumptions (\ref{h2}) and (\ref{h4}), it is worth to be noted that they are even weaker than the corresponding ones in \cite{kara2017a,kara}.
\subsubsection{The result}
\label{result-known}
Our first result states the UIB and UINN consistency of the estimators $\hat{r}_{h,\theta_0}(x)$ and $\hat{r}^\ast_{k,\theta_0}(x)$, respectively, of $r_{\theta_0}(x)$. The type of convergence considered is that of almost complete one (a.co.). Specifically, for sequences of  real random variables and positive real numbers, $\{Z_n\}$ and $\{u_n\}$, respectively, it says that $Z_n=O_{a.co.}(u_n)$ if $\exists \nu >0$ such that $\sum_{n=1}^\infty \mathbb{P}(|Z_n|>\nu u_n)<\infty$.
\begin{proposition} \label{th1}
	Let us assume that assumptions (\ref{h1}), (\ref{h2}) and (\ref{h3})-(\ref{h7}) hold.
	\begin{itemize}
		\item[(a)] If in addition Assumption (\ref{h8}) holds, then we have that
		\begin{equation*}
			\sup_{a_n\leq h\leq b_n}|\hat{r}_{h,\theta_0}(x)-r_{\theta_0}(x)|=O\left(b_n^{\beta_0}\right)+O_{a.co.}\left(\sqrt{\frac{\log n}{n\phi_{x,\theta_0}(a_n)}}\right).
		\end{equation*}
		\item[(b)] If in addition assumptions (\ref{h33bis0})-(\ref{h9}) hold, then we have that
		\begin{equation*}
			\sup_{k_{1,n}\leq k\leq k_{2,n}}|\hat{r}^\ast_{k,\theta_0}(x)-r_{\theta_0}(x)|=O\left(\phi^{-1}_{x,\theta_0}\left(\frac{k_{2,n}}{\rho_n n}\right)^{\beta_0}\right)+O_{a.co.}\left(\sqrt{\frac{\log n}{\rho_n k_{1,n}}}\right).
		\end{equation*}
	\end{itemize}
\end{proposition}

\begin{remark}
	\label{rrr}
	Proposition \ref{th1}(a) extends Theorem 3.1 in \cite{ferpv03} to the case where $h$ varies in an interval $[a_n,b_n]$ (\citealt{ferpv03} focused on the case $a_n=b_n=h$). This fact represents a very important improvement because, as it will be shown in Section \ref{data-driven}, one of the applications of Proposition \ref{th1}(a) is the validation of data-driven bandwidth selectors from an asymptotic point of view. In the same way, Proposition \ref{th1}(b) will be used in Section \ref{data-driven} to validate data-driven selectors for the number of neighbours.
\end{remark}

\subsection{The case of $\theta_0$ unknown}\label{as-2}
In practice, the direction $\theta_0$ is usually unknown and so it needs to be estimated. The results that will be presented in this section, the uniform in both bandwidth and direction (UIBD) and in both number of neighbours and direction (UINND) consistency of $\hat{r}_{h,\theta}(x)$ and $\hat{r}^\ast_{k,\theta}(x)$, respectively, play a main role to study the asymptotic behaviour of $\hat{r}_{\hat{h}}(\langle \hat{\theta}, x \rangle):=\hat{r}_{\hat{h},\hat{\theta}}(x)$ and $\hat{r}^\ast_{\hat{k}}(\langle \hat{\theta}, x \rangle):=\hat{r}^\ast_{\hat{k},\hat{\theta}}(x)$, where $\hat{h}$ and $\hat{k}$ denote some appropriate selectors for $h$ and $k$, respectively, while $\hat{\theta}$ is a suitable estimator of $\theta_0$. First, we list a few additional assumptions that we will use to state such results.

\subsubsection{Additional assumptions for UIBD and UINND consistency}
\label{addit}
\begin{description}	
	
	\item[\textit{About the space of directions.}] We assume that 
	\begin{equation}
		\textrm{card}(\Theta_n)=n^{\alpha}\quad \textrm{with} \quad \alpha>0
		\label{h11}
	\end{equation}
	and
	\begin{equation}
		\forall \ \theta \in \Theta_n, \ \left<\theta-\theta_0,\theta-\theta_0\right>^{1/2} \leq C_7 b_n.
		\label{theta-theta0}
	\end{equation}

	\item[\textit{About the functional explanatory variable.}] We assume that:
	
	\begin{equation}
		\left<X,X\right>^{1/2}\leq C_8.\label{x}
	\end{equation}

	\item[\textit{About the small-ball probability.}] We assume that:
	\begin{itemize}
		\item There exist constants $0<C_9\leq C_{10}<\infty$
		and a function $f:\mathbb{R}\longrightarrow(0,\infty)$ such that 
		\begin{equation}
			\forall \theta \in \Theta_n, \ C_9f(h)\leq\phi_{x,\theta}(h)\leq C_{10}f(h).
			\label{h33}
		\end{equation}
		(Actually, it could be the case that $f(\cdot)=f_x(\cdot)$. In the sake of brevity, we have not added the sub-index)
		
		\item There exist constants $0<C_{11}\leq C_{12}<\infty$ and sequences $\{a_n\},\{b_n\} \subset \mathbb{R}^+$ ($a_n \leq b_n$) such that, for $h \in [a_n,b_n]$ with $n$ large enough,
		\begin{equation}
			C_{11} \leq \frac{f(h/2)}{f(h)}\leq C_{12}.
			\label{h222}
		\end{equation}
		
		\item The sequences $\{a_n\}$ and $\{b_n\}$ verifies:
		\begin{equation}
			a_n\rightarrow 0 , b_n\rightarrow 0 \mbox{\text{ and }}  \frac{\log n}{n\min\left\{a_n,f(a_n)\right\}}\rightarrow 0 \label{h888}.
		\end{equation}

		\item There exist sequences $\{\rho_n\} \subset (0,1)$, $\{k_{1,n}\}\subset \mathbb{Z}^+$, $\{k_{2,n}\} \subset \mathbb{Z}^+$ ($k_{1,n}\leq k_{2,n}\leq n$) and constants $0<\lambda\leq \delta<\infty$ such that 
		\begin{equation}
			\lambda f^{-1}\left(\frac{\rho _n k_{1,n}}{n}\right)\leq\phi_{x,\theta}^{-1}\left(\frac{\rho_n k_{1,n}}{n}\right) \mbox{\text{ and }}
			\phi_{x,\theta}^{-1}\left(\frac{ k_{2,n}}{\rho_n n}\right) \leq \delta f^{-1}\left(\frac{ k_{2,n}}{\rho_n n}\right),
			\label{h33bis}
		\end{equation}
		
		\begin{equation}
			f^{-1}\left(\frac{ k_{2,n}}{\rho_n n}\right) \rightarrow 0,
			\label{h33bisbis}
		\end{equation}
		
		\begin{equation}
			\min \left\{\frac{1-\rho_n}{4}\frac{k_{1,n}}{\ln n},\frac{(1-\rho_n)^2}{4\rho_n}\frac{k_{1,n}}{\ln n}\right\} > \alpha +2 \label{hnew}
		\end{equation}
		and
		\begin{equation}
			\frac{\log n}{n \min\left\{\lambda f^{-1}(\rho_n k_{1,n}/n), f\left(\lambda f^{-1}(\rho_n k_{1,n}/n)\right) \right\}}\rightarrow 0 \label{h888bis}.
		\end{equation}

	\end{itemize}
	\item[\textit{About the kernel.}] The class of functions $\mathcal{K}_{\Theta_n}$ (see (\ref{h666})) is a pointwise measurable class such that
	\begin{equation}
		\sup_{\mathcal{Q}}\int_{0}^1\sqrt{1+\log\mathcal{N}\left(\epsilon\lVert F_{\Theta_n}\lVert_{\mathcal{Q},2},\mathcal{K}_{\Theta_n},d_{\mathcal{Q},2}\right)}d\epsilon<\infty,
		\label{h777}
	\end{equation}
	where $F_{\Theta_n}$ is the minimal envelope function of the set $\mathcal{K}_{\Theta_n}$ and the supremum is taken  over all probability measures $\mathcal{Q}$ on the measurable  space  $(\mathcal{H},\mathcal{A})$ with $||F_{\Theta_n}||_{\mathcal{Q},2}^2<\infty$.
\end{description}
Assumption (\ref{h11}) imposes that the set of directions $\Theta_n$ contains a finite number of directions, but allows it to grow up to infinity as the sample size does. In addition, taking into account the kind of results we wish to establish (UIBD and UINND consistency; see Theorem \ref{teorema theta estimado}), it is needed to impose some condition to control the bias caused by the use, in the studied statistics, of values $\theta \in \Theta_n$ other than the true value $\theta_0$. In particular, such condition should allow to link the behaviour of $d_{\theta}(\cdot,\cdot)$ and $d_{\theta_0}(\cdot,\cdot)$ (for details, see the proof of Lemma \ref{lema3}). In this paper that is done by means of Assumption (\ref{theta-theta0}). Note that, on the one hand, Assumption (\ref{theta-theta0}) implies that the larger $n$ is the closer are $\Theta_n$ and $\theta_0$; this is needed to obtain uniform consistency results on $\Theta_n$. On the other hand, the order $b_n$ in Assumption (\ref{theta-theta0}) is a technical condition (the minimal one when our proof is used) that allows to obtain the same rates of convergence as in the case of $\Theta_n=\{\theta_0\}$ (see Proposition \ref{th1}). The interested reader can find similar conditions to our Assumption (\ref{theta-theta0}) in \citet*{hardle} and \cite{xia99} (multivariate setting), and \cite{ma16} (functional setting; see also \citealt{ferr2013} for a different version of Assumption (\ref{theta-theta0})). Assumption (\ref{x}), which imposes that the explanatory variable is bounded, is not very restrictive in practice and is introduced to make the proofs easier. The role of Assumption (\ref{h33}) is to insure uniform results among all possible directions; that assumption generalizes Assumption (4) in \cite{ait} and was also used in \cite{wang2016}. Assumption (\ref{h222}) is weaker than the usual condition $0<\lim_{h\rightarrow 0}f(sh)/f(h)=\tau(s)<\infty, \ \forall s\in (0,1)$ (considered, for instance, in \citealt{kara2017a,kara}). The technical assumptions (\ref{h888}) and (\ref{h33bis})-(\ref{h888bis}) adapt those considered in \citet{kara2017a,kara} (in the context of functional nonparametric regression), respectively, to the setting of the FSIM (\ref{modelo}) (remember that, as noted in the last paragraph in Section \ref{sec-assump-known}, Assumption (H4) in \citealt{kara} should be modified in the way of our assumptions(\ref{h33bis0})-(\ref{h9})). Note that assumptions (\ref{h33})-(\ref{h888bis}) (the ones related to the small-ball probability), although technical, are not very restrictive. For instance, \cite{wang2016} showed that, under suitable conditions, $\phi_{x,\theta}(h) \approx C_{x,\theta}h$. Therefore, one can consider $f(h)=h$. Then, for such functions $\phi_{x,\theta}(\cdot)$ and $f(\cdot)$, assumptions (\ref{h33}), (\ref{h222}) and (\ref{h33bis}) are trivially verified while assumptions (\ref{h888}), (\ref{h33bisbis}) and (\ref{h888bis}) are satisfied under the conditions $\log n /(n a_n) \rightarrow 0$, $k_{2,n}/(\rho_n n) \rightarrow 0$ and $\log n /(\rho_n k_{1,n}) \rightarrow 0$, respectively. In addition, to verify Assumption (\ref{hnew}) it is sufficient that the condition $(1-\rho_n)^2 > 4(\alpha + 2) \ln n / k_{1,n}$ holds (note that none of those three conditions is restrictive and they allow the possibility that $\rho_n \rightarrow 1$). Finally, Assumption (\ref{h777}) is a natural extension of Assumption (\ref{h7}) to the current case where $\mbox{\text{card}}( \Theta_n )>1$.

\subsubsection{Main results}
Theorem \ref{teorema theta estimado} below states the UIBD and UINND consistency of $\hat{r}_{h,\theta}(x)$ and $\hat{r}^\ast_{k,\theta}(x)$, respectively, under general assumptions while, to fix the ideas, Corollary \ref{cor} shows how the rates of convergence behave in some simple case. In particular, as it will be seen along Remark \ref{nnn}, these rates of convergence are similar to the optimal ones for one-dimension  problems. This fact evidences that the main goal of constructing  procedures being insensitive to the dimensionality of the problem has been reached.
\begin{theorem}
	\label{teorema theta estimado}
	Let us assume that assumptions (\ref{h3})-(\ref{h4}), (\ref{h11})-(\ref{h222}) and (\ref{h777}) hold.
	\begin{itemize}
		\item[(a)] If in addition Assumption (\ref{h888}) holds, then we have that
		\begin{equation*}
			\sup_{\theta\in\Theta_n}\sup_{a_n\leq h\leq b_n}|\hat{r}_{h,\theta}(x)-r_{\theta_0}(x)|=O\left(b_n^{\beta_0}\right)+O_{a.co.}\left(\sqrt{\frac{\log n}{nf(a_n)}}\right).
		\end{equation*}
		
		\item[(b)] If in addition assumptions (\ref{h33bis})-(\ref{h888bis}) hold, then we have that
		\begin{equation*}
			\sup_{\theta\in\Theta_n}\sup_{k_{1,n}\leq k\leq k_{2,n}}|\hat{r}^\ast_{k,\theta}(x)-r_{\theta_0}(x)|=O\left(f^{-1}\left(\frac{k_{2,n}}{\rho_n n}\right)^{\beta_0}\right)+O_{a.co.}\left(\sqrt{\frac{\log n}{nf\left(\lambda f^{-1}(\rho_n k_{1,n}/n)\right)}}\right).
		\end{equation*}
	\end{itemize}
\end{theorem}

\begin{corollary}
	\label{cor}
	Let us assume that assumptions (\ref{h3})-(\ref{h4}), (\ref{h11})-(\ref{x}) and (\ref{h777}) hold. If in addition assumptions (\ref{h33}) and (\ref{h33bis}) hold with $f(h)=h$ and $\rho_n = \rho$ (where $0<\rho<1$ is a constant), and $k_{2,n}/n \rightarrow 0$ and $\log n /k_{1,n} \rightarrow 0$, then we have that
	\begin{equation*}
		\sup_{\theta\in\Theta_n}\sup_{k_{1,n}\leq k\leq k_{2,n}}|\hat{r}^\ast_{k,\theta}(x)-r_{\theta_0}(x)|=O\left(\left(\frac{k_{2,n}}{ n}\right)^{\beta_0}\right)+O_{a.co.}\left(\sqrt{\frac{\log n}{k_{1,n}}}\right).
	\end{equation*}
\end{corollary}

\begin{remark}
	\label{nnn}
	Theorem \ref{teorema theta estimado} extends Proposition \ref{th1} to the case where $\theta$ is unknown. As can be noted in Theorem \ref{teorema theta estimado}(b), the parameters $\lambda$ and $\rho_n$ (defined in assumptions (\ref{h33bis})-(\ref{h888bis})) affect to the rate of convergence of the $k$NN estimators. Actually, that is a consequence of having formulated our assumptions on $f(\cdot)$ in a fairly general way. Corollary \ref{cor} shows that, under the weak condition $f(h)=h$ (see the last paragraph in Section \ref{addit}), these effects  disappear. Focusing now on the specific case $f(h)=h$, let us take  $h_{o} \sim (\log n/n)^{1/(2\beta_0+1)}$, $k_{o}\sim (n^{2\beta_0}\log n)^{1/(2\beta_0+1)}$, $a_n=h_{o}-c_n$, $b_n=h_{o}+c_n$, $k_{1,n}=k_{o}-d_n$ and $k_{2,n}=k_{o}+d_n$, 
	where we have denoted $c_n= c (\log n/n)^{1/(2\beta_0+1)}$ and $d_n= c (n^{2\beta_0}\log n)^{1/(2\beta_0+1)}$, with $0<c<1$. Then one can see from Theorem \ref{teorema theta estimado}(a) and Corollary \ref{cor} that both estimates reach the rate of convergence $(\log n/n)^{\beta_0/(2\beta_0 + 1)}$ which is the well-known optimal rate for nonparametric one-dimensional problems. This attests the dimensionality reduction property of our model and estimates.
\end{remark}

\section{Data-driven parameters selection}
\label{data-driven}
An application of Theorem \ref{teorema theta estimado}(a) (Theorem \ref{teorema theta estimado}(b)) is related to the theoretical validation of both data-driven selectors for the bandwidth $h$ (for the number of neighbours $k$) and estimators for the direction $\theta_0$. Next result, which is a corollary of Theorem \ref{teorema theta estimado}, focuses on data-driven selectors based on cross-validation ideas (similar results can be derived for other usual selectors).

Let us denote $$CV(h,\theta)=n^{-1}\sum_{j=1}^n(Y_j-\hat{r}_{h,\theta}^{(-j)}(X_j))^2 \mbox{\text{ and }} CV^\ast(k,\theta)=n^{-1}\sum_{j=1}^n(Y_j-\hat{r}_{k,\theta}^{\ast(-j)}(X_j))^2,$$where, as usual, $\hat{r}_{h,\theta}^{(-j)}(\cdot)$ and $\hat{r}_{k,\theta}^{\ast(-j)}(\cdot)$ are the leave-one-out versions of $\hat{r}_{h,\theta}(\cdot)$ and $\hat{r}^\ast_{k,\theta}(\cdot)$, respectively. Then, one can considers the kernel-based estimator of $\theta_0$
$$\hat{\theta}_h=\arg \min_{\theta \in \Theta_n}CV(h,\theta)$$
(for asymptotic properties of $\hat{\theta}_h$, see \cite{ait}) and the $k$NN-based estimator
$$\hat{\theta}^\ast_k=\arg \min_{\theta \in \Theta_n}CV^\ast(k,\theta).$$ Following the same ideas, it seems natural to construct the data-driven selectors $\hat{h}$ and $\hat{k}$ as
$$\hat{h}=\arg \min_{a_n \leq h \leq b_n}CV(h,\hat{\theta}_h) \mbox{\text{ and }} \hat{k}=\arg \min_{k_{1,n} \leq k \leq k_{2,n}}CV^\ast(k,\hat{\theta}_k),$$respectively. In that way, we have two automatic estimators of $\theta_0$: one based on kernel estimation, $\hat{\theta}_{\hat{h}}$, and another based on $k$NN ideas, $\hat{\theta}^\ast_{\hat{k}}.$

\begin{corollary}
	\label{final-cor}
	\begin{itemize}
		\item[(a)] Under assumptions of Theorem \ref{teorema theta estimado}(a), we have that
		\begin{equation*}
			|\hat{r}_{\hat{h},\hat{\theta}_{\hat{h}}}(x)-r_{\theta_0}(x)|=O\left(b_n^{\beta_0}\right)+O_{a.co.}\left(\sqrt{\frac{\log n}{nf(a_n)}}\right).
		\end{equation*}
		
		\item[(b)] Under assumptions of Corollary \ref{cor}, we have that
		\begin{equation*}
			|\hat{r}^\ast_{\hat{k},\hat{\theta}^\ast_{\hat{k}}}(x)-r_{\theta_0}(x)|=O\left(\left(\frac{k_{2,n}}{ n}\right)^{\beta_0}\right)+O_{a.co.}\left(\sqrt{\frac{\log n}{k_{1,n}}}\right).
		\end{equation*}
	\end{itemize}
\end{corollary}

\begin{remark}
	\label{final-remark}
	Corollary \ref{final-cor} validates the use of cross-validation ideas to construct both estimators of the direction $\theta_0$ and data-driven selectors for the parameters $h$ and $k$ (in other words, it justifies adaptive estimation based on cross-validation ideas in the FSIM (\ref{modelo})). To the best of our knowledge, this is the first result in the literature on kernel or $k$NN adaptive estimation in the FSIM (\ref{modelo}). Actually, in the case of $k$NN adaptive estimation, there are no such kind of results not even in the multivariate single-index model.
\end{remark}

\section{Practical issues}\label{practical}
In the previous Section \ref{data-driven}, it was given theoretical validation of the estimators of the nonparametric link, $r(x)$, based on both CV-kernel and CV-$k$NN ideas, $\hat{r}_{\hat{h},\hat{\theta}_{\hat{h}}}(x)$ and $\hat{r}^\ast_{\hat{k},\hat{\theta}^\ast_{\hat{k}}}(x)$, respectively. Therefore, in practice the only additional issues that must be addressed are how to construct $\Theta_n, \ a_n, \ b_n, \ k_{1,n}$ and $k_{2,n}.$ That is the aim of this section.

\begin{description}	
	
	\item[\textit{The set of functional directions, $\Theta_n$.}] We propose to construct $\Theta_n$ in a similar way as in \cite{ait}. Specifically: 
	
	\begin{itemize}
		
		\item[(i)] Each direction $\theta \in \Theta_n$ is obtained from a $d_n$-dimensional space generated by B-spline basis functions, $\{e_1(\cdot),\ldots,e_{d_n}(\cdot)\}$. Therefore, we focus on directions
		\begin{equation}
			\theta(\cdot)=\sum_{j=1}^{d_n}\alpha_j e_j(\cdot) \ \mbox{where} \ (\alpha_1,\ldots,\alpha_{d_n}) \in \mathcal{V}. \label{theta-exp}
		\end{equation}
		\item[(ii)] The set of vectors of coefficients in (\ref{theta-exp}), $\mathcal{V}$, is obtained by means of the following procedure:
		\begin{itemize}
			
			\item[\textit{Step 1}] For each $(\beta_1,\ldots,\beta_{d_n}) \in \mathcal{C}^{d_n}$, where $\mathcal{C}=\{c_1,\ldots,c_J\} \subset \mathbb{R}^J$ denotes a set of $J$ ``seed-coefficients", construct the initial functional direction 
			\begin{equation*}
				\theta_{init}(\cdot)=\sum_{j=1}^{d_n}\beta_j e_j(\cdot). 
			\end{equation*}
			
			\item[\textit{Step 2}] For each $\theta_{init}$ in Step 1 that verifies the condition $\theta_{init}(t_0)>0$, where $t_0$ denotes a fixed value in the domain of $\theta_{init}(\cdot)$, compute $\left<\theta_{init},\theta_{init}\right>$ and construct $(\alpha_1,\ldots,\alpha_{d_n})=(\beta_1,\ldots,\beta_{d_n})/\left<\theta_{init},\theta_{init}\right>^{1/2}$.
			\item[\textit{Step 3}] Construct $\mathcal{V}$ as the set of vectors $(\alpha_1,\ldots,\alpha_{d_n})$ obtained in Step 2.
		\end{itemize}
	\end{itemize}
	Therefore, the final set of eligible functional directions is
	$$\Theta_n=\left\{\theta(\cdot)=\sum_{j=1}^{d_n}\alpha_j e_j(\cdot); \ (\alpha_1,\ldots,\alpha_{d_n}) \in \mathcal{V}\right\}.$$
	\begin{remark}
		\label{section5-remark}
		As usual, in item (i) above we consider splines of order $l\geq 1$ (degree $l-1$) and $m_n$ regularly spaced interior knots (so, $d_n=l+m_n$); note that, from the Jackson type theorem in \cite{deBoor} (page 149), if $\theta_0$ is sufficiently smooth then it is well approximated by some function in the $d_n$-dimensional space generated by B-spline basis. Note also that, by construction (see Step 2), each $\theta \in \Theta_n$ verifies the constraints $\left<\theta,\theta\right>=1$ and $\theta(t_0)>0$; so, the identifiability of the FSIM (\ref{modelo-app}) is guaranteed (for details, see Proposition 1 in \citealt{ait}). Of course, the larger $J$ in Step 1 is, the higher the size of $\Theta_n$ is (in fact, the number of initial functional directions in Step 1 is $J^{d_n}$). At this moment, it is worth being noted that our approach requires intensive computation because optimization on both $\theta$ and $h$ or $k$ is needed. Therefore, it is needed to seek for a trade-off between the size of $\Theta_n$ and the performance of the estimators. In that way, \cite{ait} suggested to consider $l=3$ and $\mathcal{C}=\{-1,0,1\}.$ 
	\end{remark}
	
	\item[\textit{The set of values for $h$: $[a_n,b_n].$}] In practice, when one needs to select some parameter (for instance $h$) via the minimization of some criterion function (for instance the CV function), it is usual to minimize over a ``wide" set, in such a way that any reasonable set of values for $h$  (for instance the set $[a_n,b_n]$ verifying the technical conditions assumed in the theoretical study) should be included in such wide set. The question of automatic selection of the interval  $[a_n,b_n]$ is still unsolved in one-dimensional nonparametric statistics, and turns in fact to be of lower importance because usually the criterion function turns to be rather flat around its minimum. Earlier references in one-dimensional setting go back to \cite{harm85} and \cite{mar85}, and the usual recommendation is to choose an interval such that the corresponding bandwidths allow to use up to $95\%$ of the sample. As we will see later along Section \ref{simulation}, this recommendation will still be efficient in the functional framework. 
	
	\item[\textit{The set of values for $k$: $\{k_{1,n}, k_{1,n}+1, \ldots, k_{2,n}\}.$}] The facts pointed just before for global kernel estimates are still valid for estimates using local bandwidths (see \citealt{vie92} for earlier advances); therefore the same recommendation can be made for the choice of the set $\{k_{1,n}, k_{1,n}+1, \ldots, k_{2,n}\}.$
\end{description}

\section{Simulation study}\label{simulation}
The aim of this section is twofold. On the one hand, to support the suggestions given in sections \ref{data-driven} and \ref{practical} related to practical issues inherent to our procedures: selection of the bandwidth ($h$) and the number of neighbours ($k$), as well as of the intervals $[a_n,b_n]$ and $[k_{1,n},k_{2,n}]$. On the other hand, to show the better performance of the $k$NN-based estimators against the one of the kernel-based estimators when heterogeneous designs are considered. 

\subsection{The design}
\label{design-sim}
For different values of $n$, observations i.i.d. $\{(X_i,Y_i)\}_{i=1}^{n+25}$ were generated from the FSIM
\begin{equation*}
	Y=r\left(\left<\theta_0,X\right>\right)+\varepsilon,
\end{equation*}
where the functional covariate was
\begin{equation}
	X(t)=a\cos(2\pi t) + b\sin(4\pi t) + 2c(t-0.25)(t-0.5) \ (t \in [0,1]).
	\label{X-sim}
\end{equation}
The same mixture distribution was considered for the random variables $a,b$ and $c$ in (\ref{X-sim}): $U(5,10)$ with probability 0.5, and $U(20,20.5)$ with probability 0.5, while each curve $X_i$ was discretized in $100$ equispaced points ($0=t_1<t_2< \cdots <t_{100}=1$). In addition, the link function was $r(u)=u^3$, the inner product was $\left\langle f, g \right\rangle=\int_{0}^{1}f(t)g(t)dt$ and $\theta_0$ was selected at random in $\Theta_n$ (more details will be given at the end of this section). Finally, $\varepsilon$ was a centred Gaussian random variable of variance equal to 0.025 times the empirical variance of $r\left(\left<\theta_0,X\right>\right)$ (i.e. signal-to-noise$=2.5\%$).

\begin{figure}[h]
	\vspace{-30pt}
	\centering
	\includegraphics[width=\textwidth]{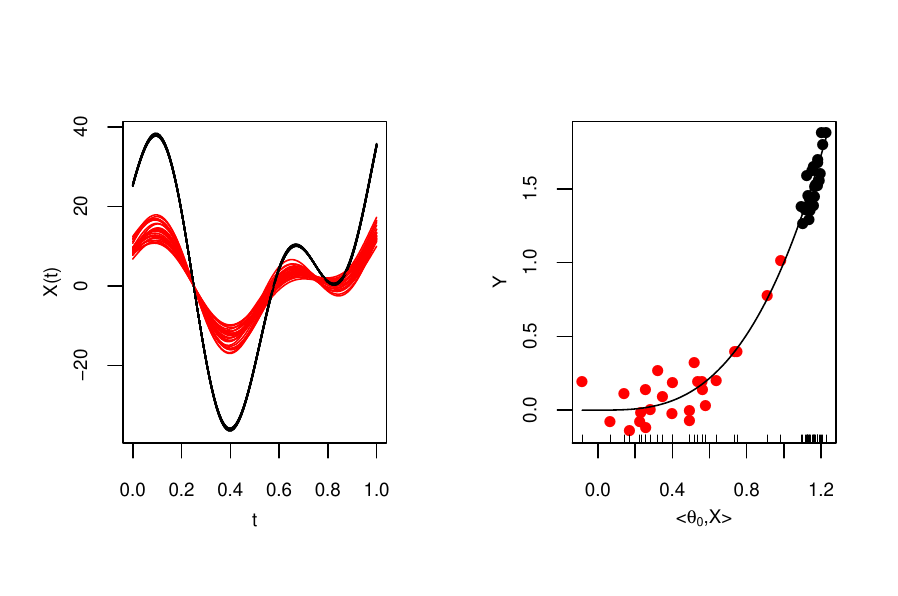}
	\vspace{-50pt}
	\caption{Sample of 50 curves $X$ (left panel) together with the corresponding scatter plot of $\{\left(\left<\theta_0,X\right>, Y\right)\}$ (right panel).}
	\label{fig1}
\end{figure}

Figure \ref{fig1} shows a sample of 50 curves (left panel) and the corresponding scatter plot of $\{\left(\left<\theta_0,X_i\right>, Y_i\right)\}_{i=1}^{50}$ (right panel). Clearly one can see two subsamples of curves, being the variability in one of them much greater than in the other. This fact gives rise to two clusters in the sample of projections, $\{\left<\theta_0,X_i\right>\}_{i=1}^{50}$; so, taking into account their location-adaptive property, one expects that the $k$NN-based estimators take advantage on the kernel-based ones.

The sample ${\cal{S}}_n=\{(X_i,Y_i)\}_{i=1}^{n+25}$ was split into two samples: a training sample, ${\cal{S}}_{n,train}=\{(X_i,Y_i)\}_{i=1}^n$, and a testing sample, ${\cal{S}}_{n,test}=\{(X_i,Y_i)\}_{i=n+1}^{n+25}$. The tuning parameters ($\hat{h}$ and $\hat{k}$) and the estimates of $\theta_0$ ($\hat{\theta}$ and $\hat{\theta}^\ast$) were constructed from the training sample by means of the cross-validation procedure proposed in Section \ref{data-driven}. The sets of functional directions ($\Theta_n$), values for $h$ ($[a_n,b_n]$) and values for $k$ ($\{k_{1,n}, k_{1,n}+1, \ldots, k_{2,n}\}$) were constructed as recommended in Section \ref{practical}. The value for $t_0$ related to $\Theta_n$ (see Step 2 in Section \ref{practical}) was fixed to $t_0=0.5$ while the considered order of the basis functions and number of interior knots were $l=3$ and $m_n=3$, respectively (as early noted, $\theta_0$ was selected at random in $\Theta_n$; once the values of $t_0$, $l$ and $m_n$ were established, we can indicate what are the coefficients of $\theta_0$: $(1.201061, 1.201061, 1.201061, 1.201061, 0,0)$; see Step 3 in Section \ref{practical}). The Epanechnikov kernel was used in the nonparametric estimates $\hat{r}(\cdot)$ and $\hat{r}^\ast(\cdot)$.

Then, the testing sample was used to measure the quality of the corresponding predictions (i.e., the performance of our procedures) through the Mean Square Error of Prediction (MSEP):
\begin{equation}
	\label{MSEP-sim}
	MSEP_n=\frac{1}{\mbox{\text{card}}({\cal{I}}_{n,test})}\sum_{i \in {\cal{I}}_{n,test}}(Y_i - \widehat{Y}_i)^2,
\end{equation}
where ${\cal{I}}_{n,test}=\{n+1,\cdots,n+25\}$ and $\widehat{Y}_i$ denotes a predicted value for $Y_i$.

\subsection{The results}
\label{results-sim}
\begin{figure}[h]
	\vspace{-30pt}
	\centering
	\includegraphics[width=\textwidth]{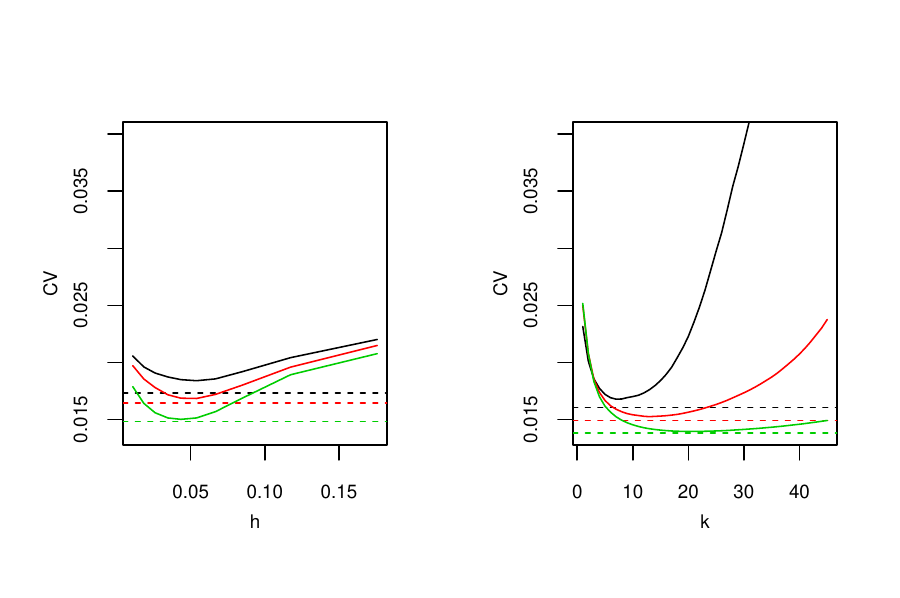}
	\vspace{-50pt}
	\caption{Average of the cross-validation functions obtained from both the kernel-based estimators (left panel) and the $k$NN-based ones as function of the bandwidth ($h$) and the number of neighbours ($k$), respectively. The dashed lines show the average of the cross-validation functions when optimal values for $h$ (left panel) and $k$ (right panel) are considered. From top to bottom, the pairs (solid curve, dashed line) correspond to $n=50,100,200$.}
	\label{fig2}
\end{figure}

For each sample size considered ($n=50,100,200$), $M=100$ replicates were generated. In order to support the suggestions given in Section \ref{practical} to construct $a_n, b_n, k_{1,n}$ and $k_{2,n}$, Figure \ref{fig2} displays the average of the cross-validation functions obtained from both the kernel-based estimator (left panel) and the $k$NN-based
estimator (right panel) when different values for the bandwidth $h$ and the number of neighbours $k$ are
considered, respectively. An interesting and practical consequence of the shown by Figure \ref{fig2} is that the optimal value for $h$ or $k$ does not change as reasonable intervals do.

\begin{figure}[h]
	\vspace{-30pt}
	\centering
	\includegraphics[width=\textwidth]{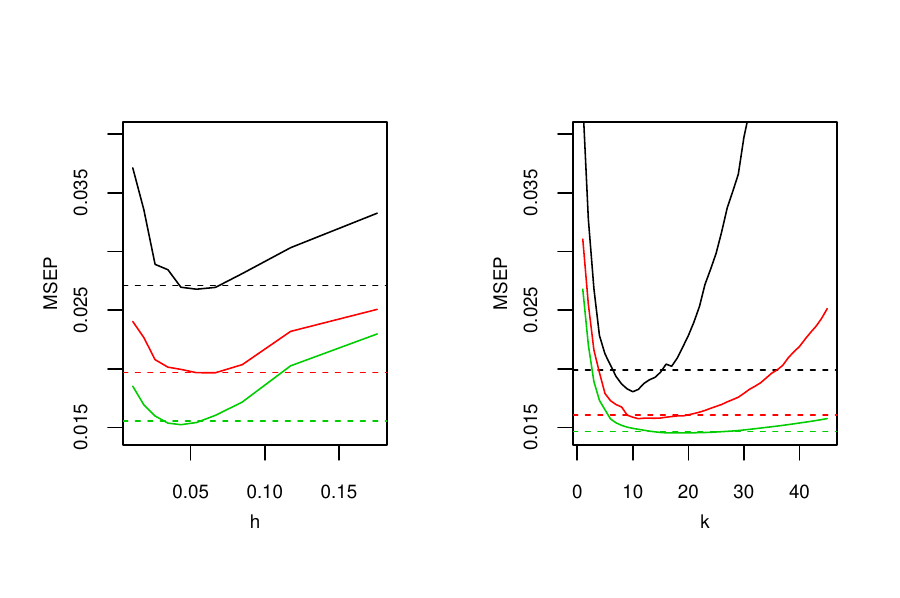}
	\vspace{-50pt}
	\caption{Average of the MSEP functions obtained from both the kernel-based estimators (left panel) and the $k$NN-based ones as function of the bandwidth ($h$) and the number of neighbours ($k$), respectively. The dashed lines show the average of the MSEP functions when values for $h$ (left panel) and $k$ (right panel) obtained from the cross-validation method are considered. From top to bottom, the pairs (solid curve, dashed line) correspond to $n=50,100,200$.}
	\label{fig3}
\end{figure}

Figure \ref{fig3} shows the average of the MSEP functions obtained from both the kernel-based estimator (left panel) and $k$NN-based
estimator (right panel) when different values for the bandwidth $h$ and the number of neighbours $k$ are
considered, respectively. The corresponding values when $h$ and $k$ are obtained from the cross-validation method are reported in Table 1. The main conclusions from Figure \ref{fig3} and Table 1 are that, for each considered sample size: (i) the
estimators are very sensitive to the values of their tuning parameters, (ii) the recommendation
given in Section \ref{practical} to construct $a_n, b_n, k_{1,n}$ and $k_{2,n}$ is appropriate (in the sense indicated in
such section), (iii) the cross-validation selectors are competitive ones, and (iv) the performance
of the $k$NN-based estimator is better than the corresponding to the kernel-based one.

\begin{center}
	\begin{tabular}{cllcclcc}
		\multicolumn{8}{l}{Table 1: Average of the MSEPs obtained when the} \\ \hline
		\multicolumn{8}{l}{CV selectors for $h$ and $k$ are used.} \\ \hline
		\multicolumn{2}{c}{$n=50$} &  & \multicolumn{2}{c}{$n=100$} &  & 
		\multicolumn{2}{c}{$n=200$} \\ \hline
		kernel & $k$NN &  & kernel & $k$NN &  & kernel & $k$NN \\ 
		$0.0271$ & $0.0199$ &  & 0.0197 & $0.0160$ &  & $0.0155$ & $0.0146$ \\ \hline
	\end{tabular}
\end{center}

\section{Application to real data}\label{application}

This section is devoted to illustrate, on a real data set, the usefulness of the FSIM (\ref{modelo}), as well as to compare the performance of the proposed adaptive kernel- and $k$NN-based estimators, $\hat{r}_{\hat{h},\hat{\theta}_{\hat{h}}}(\cdot)$ and $\hat{r}^\ast_{\hat{k},\hat{\theta}^\ast_{\hat{k}}}(\cdot)$, respectively, when the sample size increases (for details on those estimators, see Section \ref{data-driven}).

\subsection{The data} \label{data}
\begin{figure}[h]
	\centering
	\includegraphics[width=\textwidth]{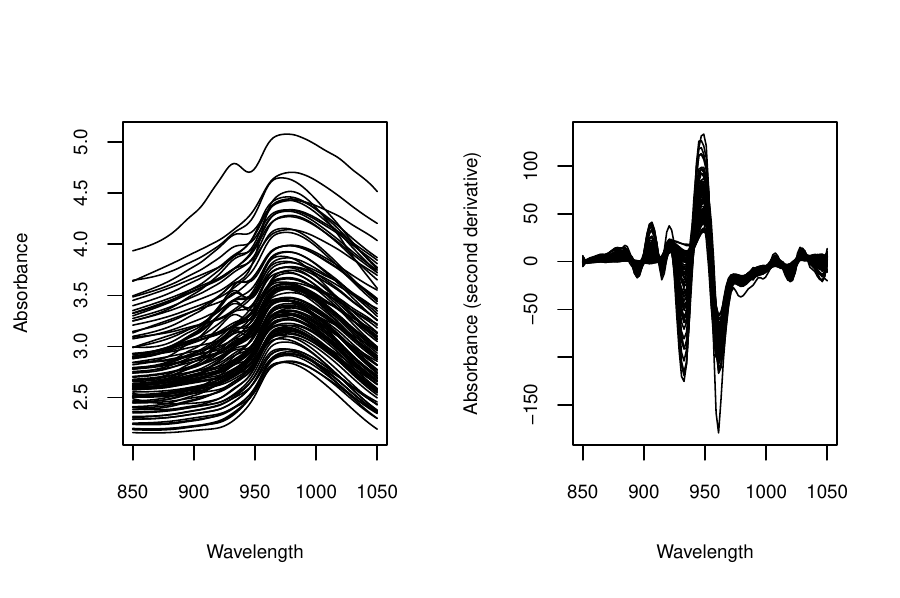}
	\caption{Sample of 100 absorbance curves $X$ (left panel) together with their second derivatives $X^{(2)}$ (right panel).}
	\label{fig4}
\end{figure}

We will analyse the well-known ``Tecator's data", a benchmark data set in the setting of nonparametric functional modelling (see, for instance, \citealt*{bur09}, \citealt{chehm11} and \citealt*{ane16} for functional pure, multiple index and sparse additive nonparametric regressions, respectively). Specifically, given 215 finely chopped pieces of meat, Tecator's data contain their corresponding fat contents ($Y_i, \ i=1,\ldots,215)$ and near-infrared absorbance spectra ($X_i, \ i=1,\ldots,215$) observed on 100 equally  wavelengths in the range $850-1050$ nm. Figure \ref{fig4} displays samples of both the absorbance curves and their second derivatives (Tecator's data are available at http://lib.stat.cmu.edu/datasets/tecator).

\subsection{The model, the tuning parameters and the data-driven selectors}
\label{consir}
As usual when one deals with Tecator's dataset, the second derivatives of the absorbance curves ($X^{(2)}$) will play the role of functional covariate. So, we focus on the FSIM
\begin{equation}
	\label{modelo-app}
	Y=r\left(\left<\theta_0,X^{(2)}\right>\right)+\varepsilon.
\end{equation}
We are interested in the performance of our procedures for different sample sizes $n$. Then, for each $n=50,100,160$, we will consider subsamples $${\cal{S}}_n=\{(X_i^{(2)},Y_i), \  i \in {\cal{I}}_n\}, \mbox{where we have denoted }  {\cal{I}}_n   =\{1,2,\ldots,n+55\}.$$ Each subsample ${\cal{S}}_n$ was split at random into two samples: a training sample, ${\cal{S}}_{n,train}=\{(X_i^{(2)},Y_i), \ i \in {\cal{I}}_{n,train}\}$, and a testing sample, ${\cal{S}}_{n,test}=\{(X_i^{(2)},Y_i), \ i \in {\cal{I}}_{n,test}\}$, where $\mbox{\text{card}}({\cal{I}}_{n,train})=n$, ${\cal{I}}_{n,train} \cup {\cal{I}}_{n,test}={\cal{I}}_n$ and ${\cal{I}}_{n,train} \cap {\cal{I}}_{n,test}=\emptyset$.

In the estimation procedures, the parameters $h$, $k$, $a_n$ , $b_n$, $k_{1,n}$ and $k_{2,n}$ were constructed from the training sample in the same way as in the simulation study (see Section \ref{practical} or Section \ref{design-sim}). Several sets of functional directions ($\Theta_n$), depending on the tuning parameter $m_n$ (number of interior knots), also were constructed as recommended in Section \ref{practical}. Values considered for $m_n$ were $2, 3, 4, 5, 6$ (note that the corresponding cardinals of $\Theta_n$ were 108, 243, 1053, 2187 and 9477, respectively), and the used value was selected by means of cross-validation ideas. The value for $t_0$ related to $\Theta_n$ (see Step 2 in Section \ref{practical}) was fixed to $t_0=(850+1050)/2$). The Epanechnikov kernel was used in the nonparametric estimates $\hat{r}(\cdot)$ and $\hat{r}^\ast(\cdot)$.

The testing sample was used to measure the quality of the corresponding predictions through the MSEP (see (\ref{MSEP-sim})).

\subsection{Performance of the procedures for different sample sizes $n$}
\label{adaptive-per}
\noindent In order to show the performance of the proposed procedures when the sample size increases, twenty partitions $({\cal{S}}_{n,train}^{(j)},{\cal{S}}_{n,test}^{(j)})$ of ${\cal{S}}_{n}$ were generated at random ($n=50,100,160; \ j=1,\ldots,20$). Then, the corresponding prediction errors, $MSEP^{(j)}_n$, were computed. Table 2 reports the average of such MSEPs.

\begin{center}
	\begin{tabular}{cllcclcc}
		\multicolumn{8}{l}{Table 2: Average of the MSEPs obtained when the} \\ \hline
		\multicolumn{8}{l}{CV selectors for $h, \ k$ and $nknots$ are used.} \\ \hline
		\multicolumn{2}{c}{$n=50$} &  & \multicolumn{2}{c}{$n=100$} &  & 
		\multicolumn{2}{c}{$n=160$} \\ \hline
		kernel & $k$NN &  & kernel & $k$NN &  & kernel & $k$NN \\ 
		$11.66$ & $10.97$ &  & 5.80 & $5.72$ &  & $4.66$ & $3.88$ \\ \hline
	\end{tabular}
\end{center}

A main suggestion from Table 2 is that, for each considered sample size, the performance of the $k$NN-based estimator is slightly better than the corresponding to the kernel-based one. In addition, the performance of each estimator improves as the sample size increases.

\subsection{The benchmark partition: Adaptive estimation in action}
\label{adaptive-action}
From now on, we focus on ${\cal{S}}_{160}$ (i.e., all the Tecator's dataset) and the partition given by ${\cal{I}}_{160,train}=\{1,2,\ldots,160\}$ and ${\cal{I}}_{160,test}=\{161,162,\ldots,215\}$. Note that this partition can be considered as a benchmark one in the sense that it is the usually considered in papers analysing the Tecator's dataset (see, for instance, \citealt{ane16}, \citealt{bur09} and \citealt{ferr2013}, among others).

\subsubsection{A Comparative study}
In a first attempt, we focus on the proposed kernel- and $k$NN-based estimates $\hat{r}(\cdot)$ and $\hat{r}^\ast(\cdot)$, respectively. In both cases, the same value for $m_n$ ($\widehat{m_n}_{CV}=4$) was selected, while the optimal bandwidth and number of neighbors where $\hat{h}_{CV}=15.80106$ and $\hat{k}_{CV}=9$, respectively. In addition,  the same estimate for $\theta_0$ ($\hat{\theta}=\hat{\theta}_{\hat{h}_{CV}}=\hat{\theta}_{\hat{k}_{CV}}^\ast$) was obtained.

\begin{figure}[h]
	\centering
	\includegraphics[width=\textwidth]{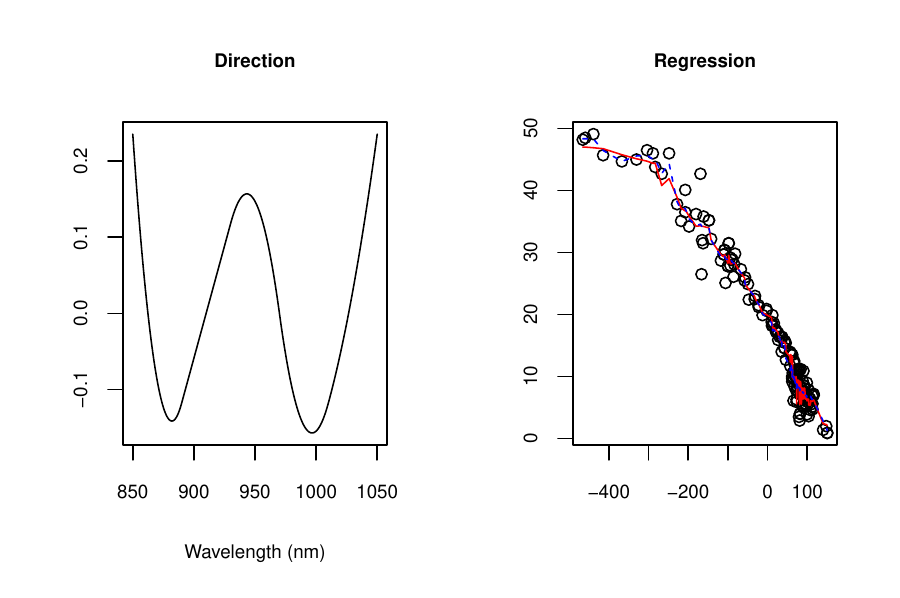}
	\caption{Left panel: Estimate of the functional direction $\theta_0$. Right panel: estimates of the regression $r(\cdot)$ by means of the $k$NN-based (solid line) and kernel-based (dashed line) estimates.}
	\label{fig5}
\end{figure}

Figure \ref{fig5} displays both the estimate for $\theta_0$ and the estimates of the regression $r(\cdot)$ (i.e., $\hat{r}_{\hat{h}_{CV}}(\cdot) $ and $\hat{r}_{\hat{k}_{CV}}^\ast(\cdot)$). On the one hand, the graphic of $\hat{\theta}$ suggests that the two bumps around wavelengths 880 and 1000, as well as the peak around wavelength 940, could be important indicators of the fat content (note that this suggestion is compatible with the findings in \citealt{ane16}). We would like to stress that one of the advantages of the FSIM against functional models dealing with the whole curves instead of with projected curves is the possibility of interpretation; as noted in the previous sentence, nice and easy interpretation is obtained in our application. On the other hand, the two estimates of the regression suggest nonlinear relationship between the fat content and the absorbance spectra (in fact, the $p$-value of the Ramsey's RESET test for linear relationship is 0.000; for details, see \citealt{rams69}). Finally, it is worth being noted the different behaviour of the considered estimates $\hat{r}_{\hat{h}_{CV}}(\cdot) $ and $\hat{r}_{\hat{k}_{CV}}^\ast(\cdot)$: in general, the kernel-based estimate is smoother than the $k$NN-based one. This fact is a consequence of two reasons: (i) the heterogeneity in the values of the covariates $\left<\hat{\theta},X_i^{(2)}\right>$, and (ii) the bandwidth ($\hat{h}_{CV}$) used in $\hat{r}_{\hat{h}_{CV}}(\cdot)$ is global (it does not depend on $x$) while the one used in $\hat{r}_{\hat{k}_{CV}}^\ast(\cdot)$ ($H_{\hat{k}_{CV},x,\hat{\theta}}$) is local (it depends on $x$). Actually, the local-adaptive bandwidth is a main appealing feature of $k$NN estimators in different settings (not only in the FSIM); in fact, as it will be shown in the remainder of this section, such feature plays a main role to attain accurate predictions.

Table 3 reports the values of the MSEPs obtained from the FSIM (\ref{modelo-app}) when estimated from both the kernel- and $k$NN-based adaptive estimators $\hat{r}_{\hat{h}_{CV}}(\cdot) $ and $\hat{r}_{\hat{k}_{CV}}^\ast(\cdot)$, respectively. The corresponding values obtained from both the functional linear model (FLM) and the pure functional nonparametric model (FNM) are also included in the table.
\begin{center}
	\begin{tabular}{llll}
		\multicolumn{4}{l}{Table 3: Values of the MSEPs from some functional models.}
		\\ \hline
		& Model & \multicolumn{2}{c}{MSEP} \\ \hline
		&  &  &  \\ 
		FLM: & $Y=\alpha _{0}+\int_{850}^{1050}X^{(2)}(t)\alpha (t)dt+\varepsilon $
		& \multicolumn{2}{c}{7.17} \\ \hline
		&  & kernel & $k$NN \\ \cline{3-4}
		FNM: & $Y=m(X^{(2)})+\varepsilon $ & $4.06$ & $1.79$ \\ 
		&  &  &  \\ 
		FSIM: & $Y=r\left( \left\langle \theta _{0},X^{(2)}\right\rangle \right)
		+\varepsilon $ & $3.49$ & $2.69$ \\ \hline
	\end{tabular}
\end{center}

In our real data application, two main conclusions can be drawn from Table 3: (i) the relationship between the fat content and the absorbance curve is nonlinear, and (ii) the FSIM estimated by means of the proposed $k$NN estimator achieves better predictive power than when it is estimated by means of the proposed kernel one. Nevertheless, the smallest value of the MSEP is obtained when the $k$NN estimator is applied to the FNM.

\subsubsection{A final boosting step}
In a second attempt, we implement a full nonparametric boosting step in the estimated FSIM. Specifically, we consider the following FNM to regress the residuals ($\hat{\varepsilon}_i$) from the FSIM on the first derivative ($X_i^{(1)}$) of the absorbance curves (the order of the derivative was selected by using cross-validation ideas):
\begin{equation}
	\hat{\varepsilon}_i=m(X_i^{(1)})+e_i, \label{second-step}
\end{equation} 
{where $e_i$ denotes the corresponding random error. Then, if $\hat{m}(\cdot)$ denotes the nonparametric estimator of $m(\cdot)$ in (\ref{second-step}), a new prediction for $Y_j$ in the test sample can be constructed as $$\hat{Y}_j=\hat{r}\left(\left<\hat{\theta},X_j^{(2)}\right>\right)+\hat{m}\left(X_j^{(1)}\right) \ (j=161,\ldots,215).$$
	
	Table 4 reports the values of the MSEP corresponding to such predictions when both functions $r(\cdot)$ and $m(\cdot)$ are estimated by means of either kernel-based or $k$NN-based estimators.
	
	\begin{center}
		\begin{tabular}{llll}
			\multicolumn{4}{l}{Table 4: Values of the MSEP when a full nonparametric} \\ 
			\multicolumn{4}{l}{boosting is applied on the residuals of the FSIM.} \\ 
			\hline
			Model & & kernel & $k$NN \\ \hline
			&& \\
			FSIM \& FNM:  & & $1.74$ & $1.53$ \\  
			
			\hline
		\end{tabular}
	\end{center}
	Several conclusions can be drawn from Table 4. On the one hand, it shows (again) the convenience of using $k$NN estimates instead of kernel ones. On the other hand, it supports the idea of considering a boosting procedure allowing to take, from the whole curve, information not captured by the functional index.
	
	\subsection{Summary of the conclusions}
	This real data analysis illustrates both the interest of the semiparametric approach and the efficiency of the $k$NN estimation procedure. On the one hand, because of its location-adaptive feature, the $k$NN approach 
	overpasses the performances of usual global smoothers, such as kernel ones, while the cross-validation procedure makes this estimate of fully automatic using. On the other hand, the semiparametric feature of the FSIM approach has the double advantage of combining interpretability of the outputs (see Figure \ref{fig5}) together with low prediction errors (see Tables 3 and 4).
	
	\section{Some tracks for future}
	\label{tracks}
	This paper has highlighted the good behaviour of the automatic and location-adaptive procedure (based on $k$NN ideas) developed in the FSIM regression setting, from both an asymptotic point of view and two finite sample size applications (a simulation study and an analysis of real data). To the best of our knowledge, this is the first work in the statistical literature in the field of data-driven location-adaptive $k$NN functional semiparametrics.
	
	Our feeling is that the ideas developed in this paper could be useful not only for the specific FSIM framework, but also one could reasonably expect that they could be used in many other settings. Keep in mind that, in multivariate analysis, the single-index model has been extended in a lot of directions. In fact, there exists a few functional literature including for instance the multiple-index model (\citealt{bouafv10}):
	\begin{equation}\label{FSIM-thetas-0}
		Y=r\left(\left<\theta_{0,1},X_1\right>+ \cdots +\left<\theta_{0,p},X_p\right>\right)+\varepsilon;
	\end{equation}
	or the functional projection pursuit regression model (\citealt{chehm11}; \citealt{ferr2013}):
	\begin{equation}\label{FSIM-thetas}
		Y=r_1\left(\left<\theta_{0,1},X\right>\right)+ \cdots +r_p\left(\left<\theta_{0,p},X\right>\right)+\varepsilon;
	\end{equation}
	or the extended model to multiple functional predictors $X^j$, as proposed  by \cite{ma16}:
	\begin{equation}
		\label{FSIM-thetas-covs-0}
		Y=r\left(\sum_{j=1}^p\left<\theta_{0,j},X^j\right>\right)+\varepsilon;
	\end{equation}
	and the functional partial linear single index model (\citealt{wang2016}):
	\begin{equation}
		\label{fplsim}
		Y=Z_1\beta_{0,1}+\cdots+Z_1\beta_{0,p}+r\left(\left<\theta_0,X\right>\right)+\varepsilon.
	\end{equation}
	The main contributions existing in the literature on the models (\ref{FSIM-thetas-0})-(\ref{fplsim}) can be summarized as follows: on the one hand, \cite{bouafv10} proved the asymptotic optimality of cross-validation-based estimators of the functional indexes $\theta_{0,j}$ in (\ref{FSIM-thetas-0}). On the other hand, \cite{chehm11}, \cite{ferr2013}, \cite{ma16} and \cite{wang2016} obtained uniform rates of convergence related to the nonparametric link functions ($r_j(\cdot)$ or $r(\cdot)$) in (\ref{FSIM-thetas}), (\ref{FSIM-thetas-covs-0}) and (\ref{fplsim}), respectively. At this moment, it is worth being noted the main difference
	between the uniform feature of our asymptotic results and those given in \cite{chehm11}, \cite{ferr2013}, \cite{ma16} and \cite{wang2016}: while the results in those papers are uniform on $x$ (or $x_j$) and $\theta$ (or $\theta_j$), our results are uniform on the tuning parameter ($h$ or $k$) and $\theta$. In addition, it is convenient to note that, while in \cite{bouafv10}, \cite{chehm11}, \cite{ferr2013} and \cite{wang2016} kernel-based estimators were studied, \cite{ma16} considered B-spline functions to approximate $\theta_{0,j}$ and $r(\cdot)$ in model (\ref{FSIM-thetas-covs-0}), and estimated both components by means of the least-squares criterion.
	
	The development of location-adaptive estimates (such as $k$NN) as well as the statement of fully automatic procedures (such as cross-validation) are, as far as we know, still missing for the general models (\ref{FSIM-thetas-0})-(\ref{fplsim}). Our guess is that the uniform ideas developed in our paper could open a way for that challenging purpose.

	\bigskip
	\noindent{\large \bf{Acknowledgements}}
	
	\noindent The authors wish to thank two anonymous referees for their helpful comments and suggestions, which greatly improved the quality of this paper. This work was supported in part by the Spanish Ministerio de Econom\'ia y Competitividad under Grant MTM2014-52876-R and
	Grant MTM2017-82724-R, in part by the Xunta de Galicia through Centro Singular de Investigaci\'on de Galicia accreditation
	under Grant ED431G/01 2016-2019 and through the Grupos de Referencia Competitiva under Grant ED431C2016-015,
	and in part by the European Union (European Regional Development Fund - ERDF). The first author also thanks the financial support from the Xunta de Galicia and the European Union (European Social Fund - ESF), the reference of which is ED481A-2018/191.

	\appendix

	\section{Appendix: Proofs}
	\label{proofs}
	From now on, $C$ denotes a generic positive constant which may take different values from one formula to another.
	
	Before presenting the proofs of our Proposition \ref{th1}, Theorem \ref{teorema theta estimado}, Corollary \ref{cor} and Corollary \ref{final-cor}, we first enunciate some known auxiliary results that play a main role in our proofs. In such results, $Z_1, Z_2, \dots, Z_n$ are i.i.d. variables taking values in a measurable space $(\mathcal{Z},\mathcal{A})$ and $\mathcal{K}$ is a pointwise measurable class of functions $\{g:\mathcal{Z}\longrightarrow\mathbb{R\}}$ with envelope function $F$. In addition, we denote
	
	\[
	\alpha_n(g)=\frac{1}{\sqrt{n}}\sum_{i=1}^n\left(g(Z_i)-\mathbb{E}(g(Z_i)\right)), \ ||\alpha_n(g)||_{\mathcal{K}}=\sup_{g \in \mathcal{K}}|\alpha_n(g)|, \ \lvert\lvert\cdot\rvert\rvert_p=\sqrt[p]{\mathbb{E}(\cdot)^p}
	\]	
	and
	\[
	J(1,\mathcal{K})= \sup_{\mathcal{Q}}\int_{0}^1\sqrt{1+\log\mathcal{N}\left(\epsilon||F||_{\mathcal{Q},2},\mathcal{K},d_{\mathcal{Q},2}\right)}d\epsilon,
	\]
	where the supremum is taken over all the probability measures $\mathcal{Q}$ on the measure space $(\mathcal{Z},\mathcal{A})$ with $||F||_{\mathcal{Q},2}<\infty$. (For additional notation, see Section \ref{asymptotics})

	\subsection{Some auxiliary results}
	
	\begin{lemma}[\bf Theorem 2.14.1 in \citealt*{van1196}, p. 239] \label{lemaA.1}	We have that:
		\begin{equation*}
			\left\lvert\left\lvert \hspace{1mm} \left\lvert\left\lvert \alpha_n(g) \right\rvert \right\rvert_{\mathcal{K}} \right\rvert\right\rvert_p\leq C J(1,\mathcal{K}) \left\lvert \left\lvert F\right \rvert\right \rvert_{p\vee2},
		\end{equation*}
		where $s\vee t$ is the spermium of $s$ and $t$.
		
	\end{lemma}

	\begin{lemma}[\bf Theorem 3.1 in \citealt*{dony}, p. 314] \label{lemaA.2} If the class $\mathcal{K}$ is such that $\mathbb{E}\left\lvert\left\lvert \alpha_n(g) \right\rvert \right\rvert_{\mathcal{K}}\leq C\left\lvert \left\lvert F\right \rvert\right \rvert_2$, then, for any $A\in\mathcal{A}$, we have:
		\begin{equation*}
			\mathbb{E}\left\lvert\left\lvert \alpha_n(g1_A) \right\rvert \right\rvert_{\mathcal{K}} \leq 2C \left\lvert \left\lvert F1_A\right \rvert\right \rvert_2.
		\end{equation*}
	\end{lemma}

	\begin{lemma}[\bf Bernstein type inequality in \citealt*{dony}, p. 321] \label{lemaA.3}	Assume that the variables $Z_1,Z_2,\dots,Z_n$ satisfy for some $H>0$,
		\begin{equation*}
			\mathbb{E}\left(F^p(Z)\right)\leq\frac{p!}{2}\sigma^2H^{p-2},
		\end{equation*}
		where $\sigma^2\geq\mathbb{E}(F^2(Z))$. Then, by denoting $\beta_n=\mathbb{E}(\lvert\lvert\sqrt{n}\alpha_n(g)\rvert\rvert_{\mathcal{K}})$ we have for any $t>0$:
		\begin{equation*}
			\mathbb{P}\left\{ \max_{1\leq k\leq n} \left\lvert\left\lvert \sqrt{k}\alpha_k(g) \right\rvert \right\rvert_{\mathcal{K}} \geq \beta_n+t \right\}\leq \exp\left(-\frac{t^2}{2n\sigma^2+2tH}\right).
		\end{equation*}	
	\end{lemma}

	\begin{lemma}[\bf Lema 6.1 in \citealt{kara}, p. 186] \label{lemaA.4}
		Let $X_1,\dots,X_n$ be independent Bernoulli random variables with $\mathbb{P}(X_i)=p$ for all $i=1,\dots,n$. Set $U=X_1+\dots+X_n$ and $\mu=pn$. Then, for any $w>0$, we have:
		\begin{equation*}
			\mathbb{P}\left(U\geq (1+w)\mu\right)\leq \exp\{-\mu\min\{w,w^2\}/4\},
		\end{equation*}
		and if $w\in(0,1)$, we have
		\begin{equation*}
			\mathbb{P}\left(U\leq (1-w)\mu\right)\leq \exp\{-\mu w^2/2\}.
		\end{equation*}
	\end{lemma}
	
	\subsection{Proof of Proposition \ref{th1}}
	
	Results in Proposition \ref{th1}(a) and Proposition \ref{th1}(b) are direct consequence of Theorems 3.1 in \cite{kara2017a,kara}, respectively. On the one hand, one must note that, when
	$\theta_0$ is known, $\hat{r}_{h,\theta_0}(\cdot)$ and $\hat{r}^\ast_{k,\theta_0}(\cdot)$ are
	kernel- and \textit{k}NN-type estimators, respectively, based on the semi-metric $d_{\theta_0}(\cdot,\cdot)$, of the nonparametric regression operator, $r_{\theta_0}(\cdot)$, between the scalar variable $Y$ and the functional covariate $X$. 	On the one hand, in the case of Proposition \ref{th1}(b), one must take into account the correction relative to the rate of convergence in Theorem 3.1 in \cite{kara} indicated in the last paragrhapn in Section \ref{sec-assump-known}.
	
	Actually, our assumptions (\ref{h2}) and (\ref{h4}) are slightly different (weaker) of assumptions (6) and (10) in \cite{kara} and assumptions H1 and H3 in \cite{kara2017a}; to show that their Theorems 3.1 maintain when one uses our assumptions instead of the corresponding ones in \cite{kara2017a,kara}, it is sufficient to prove Corollary 3.3 in \cite{kara2017a} following the proof of our Corollary \ref{lema2} (see below).
	
	\subsection{Proof of Theorem \ref{teorema theta estimado}(a)}
	
	We will follow the scheme used in \cite{kara2017a}, who focused on the UIB consistency of the kernel estimator of the nonparametric regression of a scalar response on a functional explanatory variable. Although our theorem differs, with respect to the one in \cite{kara2017a}, in both the model and the type of consistency to prove (we focus on the FSIM (\ref{modelo}) instead of the functional nonparametric model and our aim is the UIBD consistency instead of the UIB one), their scheme of proof can be followed once one adapts the assumptions in a suitable way.
	
	Taking into account that
	\begin{equation*}
		\hat{r}_{h,\theta}(x)=\frac{\hat{g}_{h,\theta}(x)}{\hat{F}_{h,\theta}(x)},
	\end{equation*}	
	where we have denoted
	\begin{equation*}
		\hat{g}_{h,\theta}(x)=\frac{1}{n\phi_{x,\theta}(h)}\sum_{i=1}^nK\left(h^{-1}d_{\theta}\left(x,X_i\right)\right)Y_i \mbox{\text{ and }} \hat{F}_{h,\theta}(x)=\frac{1}{n\phi_{x,\theta}(h)}\sum_{i=1}^nK\left(h^{-1}d_{\theta}\left(x,X_i\right)\right),
	\end{equation*}
	we can write
	\begin{equation*}
		\hat{r}_{h,\theta}(x)-r_{\theta_0}(x)=\hat{B}_{h,\theta}(x)+\frac{\hat{R}_{h,\theta}(x)}{\hat{F}_{h,\theta}(x)}+\frac{\hat{Q}_{h,\theta}(x)}{\hat{F}_{h,\theta}(x)},
	\end{equation*}	
	where
	\begin{equation*}
		\hat{B}_{h,\theta}(x)=\frac{\mathbb{E}\left(\hat{g}_{h,\theta}(x)\right)}{\mathbb{E}\left(\hat{F}_{h,\theta}(x)\right)}-r_{\theta_0}(x), \ \hat{R}_{h,\theta}(x)=-\hat{B}_{h,\theta}(x)\left(\hat{F}_{h,\theta}(x)-\mathbb{E}\left(\hat{F}_{h,\theta}(x)\right)\right)
	\end{equation*}	
	and
	\begin{equation*}
		\hat{Q}_{h,\theta}(x)=\left(\hat{g}_{h,\theta}(x)-\mathbb{E}\left(\hat{g}_{h,\theta}(x)\right)\right)-r_{\theta_0}(x)\left(\hat{F}_{h,\theta}(x)-\mathbb{E}\left(\hat{F}_{h,\theta}(x)\right)\right).
	\end{equation*}
	Thus, the proof of our Theorem \ref{teorema theta estimado}(a) is completed once we prove the following four results.
	
	\begin{lemma} Under assumptions (\ref{h4}), (\ref{h11}), (\ref{h33}), (\ref{h888}) and (\ref{h777}), we have that:
		\label{lema1}
		\begin{equation*}
			\sup_{\theta\in\Theta_n}\sup_{a_n\leq h\leq b_n}\left|\hat{F}_{h,\theta}(x)-\mathbb{E}\left(\hat{F}_{h,\theta}(x)\right)\right|=O_{a.co.}\left(\sqrt{\frac{\log n}{nf(a_n)}}\right).
		\end{equation*}
	\end{lemma}
	
	\begin{corollary} Under assumptions of Lemma \ref{lema1} together with Assumption (\ref{h222}), there exists $C>0$ such that
		\label{lema2}
		\begin{equation*}
			\sum_{n=1}^{\infty}	\mathbb{P}\left(\inf_{\theta\in\Theta_n}\inf_{a_n\leq h\leq b_n}\hat{F}_{h,\theta}(x)<C\right)<\infty.
		\end{equation*}
	\end{corollary}
	
	\begin{lemma} Under assumptions and (\ref{h3}), (\ref{h4}) and (\ref{theta-theta0})--(\ref{h222}), we have that:
		\label{lema3}
		\begin{equation*}
			\sup_{\theta\in\Theta_n}\sup_{a_n\leq h\leq b_n}\left|\hat{B}_{h,\theta}(x)\right|=O\left(b_n^{\beta_0}\right).
		\end{equation*}
	\end{lemma}

	\begin{lemma} Under assumptions (\ref{h}), (\ref{h4}), (\ref{h11}), (\ref{h33}), (\ref{h888}) and (\ref{h777}), we have that:
		\label{lema4}
		\begin{equation*}
			\sup_{\theta\in\Theta_n}\sup_{a_n\leq h\leq b_n}\left|\hat{g}_{h,\theta}(x)-\mathbb{E}\left(\hat{g}_{h,\theta}(x)\right)\right|=O_{a.co.}\left(\sqrt{\frac{\log n}{nf(a_n)}}\right).
		\end{equation*}
	\end{lemma}

	\subsubsection{ Proof of Lemma \ref{lema1}}
	Taking Assumption (\ref{h33}) into account, it suffices to prove that there exist $\eta_0> 0$ and $b_0>0$ such that	
	\begin{equation}
		\sum_{n=1}^{\infty}\mathbb{P}\left(\sup_{\theta\in\Theta_n}\sup_{a_n\leq h\leq b_0}\sqrt{\frac{n\phi_{x,\theta}(a_n)}{\log n}}\left|\hat{F}_{h,\theta}(x)-\mathbb{E}\left(\hat{F}_{h,\theta}(x)\right)\right|\geq \eta_0 \right)<\infty. \label{a0a}
	\end{equation}
	Let us define
	\begin{equation*}
		h_j=2^ja_n,\ L(n)=\max\{j:h_j\leq 2b_0\} \textrm{ and } K_{h,\theta}(X_i)=K(h^{-1}d_\theta(x,X_i).
	\end{equation*}
	In addition, let us denote
	$$\Theta_n=\{\theta_1,\dots,\theta_{n^{\alpha}}\}, \ \alpha_{n}(g)=\frac{1}{\sqrt{n}}\sum_{i=1}^n\left(g(X_i)-\mathbb{E}\left(g(X_i)\right)\right)$$
	and, for $1\leq j \leq L(n)$ and $1\leq m \leq n^\alpha$,
	\begin{equation*}
		\mathcal{G}_{j,m}=\left\{\cdot\longrightarrow K\left(h^{-1}d_{\theta_m}(x,\cdot)\right) \textrm{ where } h_{j-1}\leq h\leq h_j\right\}.
	\end{equation*}
	Thus, taking into account that
	\begin{equation*}
		\hat{F}_{h,\theta}(x)-\mathbb{E}\left(\hat{F}_{h,\theta}(x)\right)=\frac{1}{\sqrt{n}\phi_{x,\theta}(h)}\alpha_{n}(K_{h,\theta}),
	\end{equation*}
	we can write
	\begin{eqnarray}
		&&\mathbb{P}\left(\sup_{\theta\in\Theta_n}\sup_{a_n\leq h\leq b_0}\sqrt{\frac{n\phi_{x,\theta}(a_n)}{\log n}}\left|\hat{F}_{h,\theta}(x)-\mathbb{E}\left(\hat{F}_{h,\theta}(x)\right)\right|\geq \eta_0 \right)\nonumber\\
		&\leq&\sum_{m=1}^{n^{\alpha}}\sum_{j=1}^{L(n)}\mathbb{P}\left(\frac{1}{\sqrt{n\log n\phi_{x,\theta_m}\left(h_j/2\right)}}\left\lVert\sqrt{n}\alpha_{n}(g)\right\rVert_{\mathcal{G}_{j,m}}\geq \eta_0 \right)\nonumber\\
		&\leq& n^{\alpha}L(n)\max_{m=1,\dots,n^{\alpha}}\max_{j=1,\dots, L(n)}\mathbb{P}\left(\max_{1\leq k\leq n}\left\lVert\sqrt{k}\alpha_{k}(g)\right\rVert_{\mathcal{G}_{j,m}}\geq \eta_0\sqrt{n\log n\phi_{x,\theta_m}\left(h_j/2\right)} \right).\label{aa}
	\end{eqnarray}
	In order to bound the probability that appears in (\ref{aa}) by means of the the Bernstein's inequality formulated in Lemma \ref{lemaA.3}, we first study the asymptotic behaviour of
	\begin{equation*}
		\sigma^2=\mathbb{E}\left(G_{j,m}^2(X)\right) \textrm{ and } \beta_n=\mathbb{E}\left(\left\lVert\sqrt{n}\alpha_{n}(g)\right\rVert_{\mathcal{G}_{j,m}}\right),
	\end{equation*}
	where $G_{j,m}$ denotes the minimal envelope function of the class $\mathcal{G}_{j,m}$. It follows from Assumption (\ref{h4}) that
	\begin{equation*}
		G_{j,m}(z)\leq C_6 1_{B_{\theta_m}(x,h_j/2)}(z).
	\end{equation*}
	
	Hence, for all $p\geq 1$, we have that
	\begin{equation*}
		\mathbb{E}\left(G_{j,m}(X)^p\right)\leq C_6^p \phi_{x,\theta_m}\left(h_j/2\right).
	\end{equation*}
	In particular,
	\begin{equation*}
		\sigma^2=O\left( \phi_{x,\theta_m}\left(h_j/2\right)\right)
	\end{equation*}
	holds. Focusing now on $\beta_n$, we obtain, by combining Assumption (\ref{h777}) together with Lemma \ref{lemaA.1}, that
	\begin{equation*}
		\mathbb{E}\left(\left\lVert\alpha_{n}(g)\right\rVert_{\mathcal{G}_{j,m}}\right) \leq \mathbb{E}\left(\left\lVert\alpha_{n}(g)\right\rVert_{\mathcal{K}_{\Theta_n}}\right)\leq CJ(1,\mathcal{K}_{\Theta_n})\left\lVert F_{\Theta _ n} \right\rVert_2\leq C\left\lVert F_{\Theta _ n} \right\rVert_2.
	\end{equation*}
	Thus, the conditions of  Lemma \ref{lemaA.2} are verified for the class $\mathcal{G}_{j,m}$
	and the envelope function $F_{\Theta _ n}$ (note that, in particular, $F_{\Theta _ n}$ is an envelope function of the class $\mathcal{G}_{j,m}$). So, from such lemma it follows that
	\begin{equation*}
		\mathbb{E}\left(\left\lVert\alpha_{n}\left(g1_{B_{\theta_m}\left(x,h_j/2\right)}\right)\right\rVert_{\mathcal{G}_{j,m}}\right)\leq C\left\lVert F1_{B_{\theta_m}\left(x, h_j/2\right)}\right\rVert_2.
	\end{equation*}
	Finally, taking into account (\ref{h4}), we obtain that:
	\begin{equation*}
		\beta_n=\mathbb{E}\left(\left\lVert\sqrt{n}\alpha_{n}(g)\right\rVert_{\mathcal{G}_{j,m}}\right)=\mathbb{E}\left(\left\lVert\sqrt{n}\alpha_{n}\left(g1_{B_{\theta_m}\left(x,h_j/2\right)}\right)\right\rVert_{\mathcal{G}_{j,m}}\right)
		\leq C\sqrt{n\phi_{x,\theta_m}\left(h_j/2\right)}.
	\end{equation*}
	Now, we can apply the Bernstein's inequality (see Lemma \ref{lemaA.3}) with
	\begin{equation*}
		\beta_n=O\left(\sqrt{n\phi_{x,\theta_m}\left(h_j/2\right)}\right), \ \sigma^2=O\left( \phi_{x,\theta_m}\left(h_j/2\right)\right) \mbox{\text{ and }}  t=\frac{\eta_0}{2} \sqrt{n \log n\phi_{x,\theta_m}\left(h_j/2\right)},
	\end{equation*}
	and we obtain:
	\begin{equation}
		\begin{array}{c}
			\mathbb{P}\left(\underset{1\leq k\leq n}{\max}\left\lVert\sqrt{k}\alpha_{k}(g)\right\rVert_{\mathcal{G}_{j,m}}\geq \eta_0\sqrt{n\log n\phi_{x,\theta_m}\left(h_j/2\right)} \right) \leq \mathbb{P}\left(\underset{1\leq k\leq n}{\max}\left\lVert\sqrt{k}\alpha_{k}(g)\right\rVert_{\mathcal{G}_{j,m}}\geq \beta_n+t \right) \\ 
			\leq \exp\left\{-\eta_0^2\frac{\log n}{8+C\sqrt{\frac{\log n}{n\phi_{x,\theta_m}(h_j/2)}}}\right\}
			\leq \exp\left\{-\eta_0^2\frac{\log n}{8+C\sqrt{\frac{\log n}{nf(h_j/2)}}}\right\}
			\leq n^{-C'\eta_0^2}
		\end{array}
		\label{bb}
	\end{equation}
	(Note that the penultimate inequality is consequence of Assumption (\ref{h33}) and the last one of Assumption (\ref{h888})) Moreover, from (\ref{aa}) and (\ref{bb}) together with the fact that $L(n)\leq2\log n$, we get that
	\begin{equation}
		\mathbb{P}\left(\sup_{\theta\in\Theta_n}\sup_{a_n\leq h\leq b_0}\sqrt{\frac{n\phi_{x,\theta}(a_n)}{\log n}}\left|\hat{F}_{h,\theta}(x)-\mathbb{E}\left(\hat{F}_{h,\theta}(x)\right)\right|\geq \eta_0 \right)\leq C n^{-C'\eta_0^2+\alpha}\log n .\label{b0b}
	\end{equation}
	Finally, by choosing now $\eta_0$ such that $C'\eta_0^2-\alpha>1$, (\ref{a0a}) follows from (\ref{b0b}).
	
	\subsubsection{Proof of Corollary \ref{lema2}}
	From assumptions (\ref{h4}), (\ref{h33}) and (\ref{h222}), we obtain that, for $n$ large enough,
	\begin{equation}
		\mathbb{E}\left(\hat{F}_{h,\theta}(x)\right)\geq C_5 \frac{\phi_{x,\theta}(h/2)}{\phi_{x,\theta}(h)}\geq  \frac{C_5C_9}{C_{10}}\frac{f(h/2)}{f(h)} \geq \frac{C_5C_9C_{11}}{C_{10}}=C'>0, \ \forall h\in [a_n,b_n] \mbox{\text{ and }} \forall \theta \in \Theta_n. \label{EF}
	\end{equation}
	Thus, denoting $C=C'/2$, it verifies that
	\begin{equation*}
		\mathbb{P}\left(\inf_{\theta\in\Theta_n}\inf_{h\in [a_n,b_n]}\hat{F}_{h,\theta}(x)\leq C\right)\leq \mathbb{P}\left(\sup_{\theta\in\Theta_n}\sup_{h\in [a_n,b_n]}\left|\mathbb{E}\left(\hat{F}_{h,\theta}(x)\right)-\hat{F}_{h,\theta}(x)\right|\geq C\right),
	\end{equation*}
	and Lemma \ref{lema1} leads to the desired result.\\
	
	
	\subsubsection{ Proof of Lemma \ref{lema3}}
	We have that
	\begin{eqnarray}
		|\hat{B}_{h,\theta}(x)\mathbb{E}(\hat{F}_{h,\theta}(x))|
		&=&|\frac{1}{\phi_{x,\theta}(h)}\mathbb{E}\left[ K(h^{-1}d_{\theta}(x,X))\left(\mathbb{E}(Y|X)-r_{\theta_0}(x)\right)\right]|\nonumber\\
		&=&|\frac{1}{\phi_{x,\theta}(h)}\mathbb{E}\left[ K(h^{-1}d_{\theta}(x,X))\left(r_{\theta_0}(X)-r_{\theta_0}(x)\right)\right]|\nonumber\\
		&\leq& C_6\frac{1}{\phi_{x,\theta}(h)}\mathbb{E}\left[1_{B_{\theta}\left(x,h/2\right)}(X) d_{\theta_0}(X,x)^{\beta_0}\right]. \label{ine}
	\end{eqnarray}
	(Note that the inequality in (\ref{ine}) is consequence of assumptions (\ref{h3}) and (\ref{h4})) In addition, assumptions (\ref{theta-theta0}) and (\ref{x}) allow us to write that, if $d_{\theta}(X,x)<h/2$ holds, then, for all $h \in [a_n,b_n]$,
	\begin{eqnarray}
		d_{\theta_0}(X,x)&=& d_{\theta_0}(X,x)-d_{\theta}(X,x)+d_{\theta}(X,x) \leq \left|\left<X-x,\theta_0-\theta\right>\right|+d_{\theta}(X,x)\nonumber\\
		&\leq & \left<X-x, X-x\right>^{1/2}  \left<\theta_0-\theta,\theta_0-\theta\right>^{1/2}+d_{\theta}(X,x) \leq 2C_7 C_8 b_n+h/2 \leq Cb_n. \label{dist0}
	\end{eqnarray}
	Now, from (\ref{ine}) and (\ref{dist0}) together with assumptions (\ref{h33}) and (\ref{h222}), we obtain that, for all $h \in [a_n, b_n]$ and for all $\theta \in \Theta_n$,
	\begin{equation}
		|\hat{B}_{h,\theta}(x)\mathbb{E}(\hat{F}_{h,\theta}(x))| \leq C \frac{\phi_{x,\theta}(h/2)}{\phi_{x,\theta}(h)}b_n^{\beta_0} \leq C \frac{f(h/2)}{f(h)}b_n^{\beta_0} \leq Cb_n^{\beta_0}. \label{BF}
	\end{equation}
	Finally, (\ref{EF}) and (\ref{BF}) complete the proof.

	\subsubsection{ Proof of Lemma \ref{lema4}}
	Taking Assumption (\ref{h33}) into account, it suffices to prove that there exist $\eta'_0> 0$ and $b_0>0$ such that	
	\begin{equation}
		\sum_{n=1}^{\infty}\mathbb{P}\left\{\sup_{\theta\in\Theta_n}\sup_{a_n\leq h\leq b_0}\sqrt{\frac{n\phi_{x,\theta}(a_n)}{\log n}}\left|\hat{g}_{h,\theta}(x)-\mathbb{E}\left(\hat{g}_{h,\theta}(x)\right)\right|\geq \eta_0' \right\}<\infty. \label{ppp}
	\end{equation}
	As in the proof of Lemma \ref{lema1}, let us define $h_j=2^ja_n \textrm{ and } L(n)=\max\{j:h_j\leq 2b_0\}$. In addition, let us denote
	$$\alpha'_{n}(g)=\frac{1}{\sqrt{n}}\sum_{i=1}^n\left(Y_ig(X_i)-\mathbb{E}\left(Y_ig(X_i)\right)\right), \ \beta_n'=\mathbb{E}\left(\left\lVert\sqrt{n}\alpha_{n}'(g)\right\rVert_{\mathcal{G}_{j,m}'}\right) \mbox{\text{ and }} \sigma'^2= \mathbb{E}\left(G_{j,m}'(X,Y)^2\right),$$
	where
	\begin{equation*}
		\mathcal{G}_{j,m}'=\left\{(z,y)\longrightarrow yK\left(h^{-1}d_{\theta_m}(x,z)\right) \textrm{ where } h_{j-1}\leq h\leq h_j\right\}
	\end{equation*}
	and $G_{j,m}'$ denotes the minimal envelope function of the class $\mathcal{G}_{j,m}'$ ($1\leq j \leq L(n)$ and $1\leq m \leq n^\alpha$).
	
	In a similar way as in the proof of Lemma \ref{lema1}, we have that
	\begin{eqnarray}
		&&\mathbb{P}\left\{\sup_{\theta\in\Theta_n}\sup_{a_n\leq h\leq b_0}\sqrt{\frac{n\phi_{x,\theta}(a_n)}{\log n}}\left|\hat{g}_{h,\theta}(x)-\mathbb{E}\left(\hat{g}_{h,\theta}(x)\right)\right|\geq \eta_0' \right\}\nonumber\\
		&\leq& n^{\alpha}L(n)\max_{m=1,\dots,n^{\alpha}}\max_{j=1,\dots, L(n)}\mathbb{P}\left(\max_{1\leq k\leq n}\left\lVert\sqrt{k}\alpha'_{k}(g)\right\rVert_{\mathcal{G}_{j,m}}\geq \eta'_0\sqrt{n\log n\phi_{x,\theta_m}\left(h_j/2\right)} \right).\label{uu}
	\end{eqnarray}
	In addition, by using the same ideas as for the proof of Lemma \ref{lema1}, and taking into account the Assumption (\ref{h}), we get:
	\begin{equation*}
		\mathbb{E}\left(G_{j,m}'(X,Y)^p\right)\leq C^p \phi_{x,\theta_m}\left(h_j/2\right), \ \sigma'^2=O\left( \phi_{x,\theta_m}\left(h_j/2\right)\right)
		\mbox{\text{ and }} \beta_n'=O\left(\sqrt{n\phi_{x,\theta_m}\left(h_j/2\right)}\right).
	\end{equation*}
	Now, from the Bernstein's inequality (see Lemma \ref{lemaA.3}), we obtain
	\begin{eqnarray}
		\mathbb{P}\left\{\max_{1\leq k\leq n}\left\lVert\sqrt{k}\alpha_{k}'(g)\right\rVert_{\mathcal{G}_{j,m}'}\geq \eta_0'\sqrt{n\log n\phi_{x,\theta_m}\left(h_j/2\right)} \right\}&\leq& \mathbb{P}\left\{\max_{1\leq k\leq n}\left\lVert\sqrt{k}\alpha_{k}'(g)\right\rVert_{\mathcal{G}_{j,m}'}\geq \beta_n'+t \right\}\nonumber\\
		&\leq&n^{-C'\eta_0'^2},\label{ttt}
	\end{eqnarray}
	while from (\ref{uu}) and (\ref{ttt}) we have that
	\begin{equation}
		\mathbb{P}\left\{\sup_{\theta\in\Theta_n}\sup_{a_n\leq h\leq b_0}\sqrt{\frac{n\phi_{x,\theta}(a_n)}{\log n}}\left|\hat{g}_{h,\theta}(x)-\mathbb{E}\left(\hat{g}_{h,\theta}(x)\right)\right|\geq \eta_0' \right\} \leq n^{\alpha} n^{-C'\eta_0'^2}\log n. \label{fin}
	\end{equation}
	Finally, by choosing $\eta_0'$ such that $C'\eta_0'^2-\alpha>1$, (\ref{ppp}) follows from (\ref{fin}).
	
	\subsection{Proof of Theorem \ref{teorema theta estimado}(b)}
	
	We will follow the scheme of \cite{kara}, but taking into account that in our setting, for a fixed $k$, the random bandwidth also depends on $\theta$.
	
	We have that
	\begin{eqnarray}
		\sup_{\theta\in \Theta_n}\sup_{k_{1,n}\leq k\leq k_{2,n}}\left|\hat{r}^\ast_{k,\theta}(x)-r_{\theta_0}(x)\right|&=&
		\sup_{\theta\in \Theta_n}\sup_{k_{1,n}\leq k\leq k_{2,n}}\left|\hat{r}^\ast_{k,\theta}(x)-r_{\theta_0}(x)\right|1_{\left\{\phi_{x,\theta}^{-1}\left(\frac{\rho_n k_{1,n}}{n}\right)\leq H_{k,x,\theta}\leq \phi_{x,\theta}^{-1}\left(\frac{k_{2,n}}{\rho_n n}\right)\right\}}\nonumber \\ 
		&+&\sup_{\theta\in \Theta_n}\sup_{k_{1,n}\leq k\leq k_{2,n}}\left|\hat{r}^\ast_{k,\theta}(x)-r_{\theta_0}(x)\right|1_{\left\{H_{k,x,\theta}\not\in\left(\phi_{x,\theta}^{-1}\left(\frac{\rho_n k_{1,n}}{n}\right), \phi_{x,\theta}^{-1}\left(\frac{ k_{2,n}}{\rho_n n}\right)\right)\right\}},\nonumber
	\end{eqnarray}	
	where $\rho_n \in (0,1)$ was defined in Assumption (\ref{h33bis}). Thus, taking Assumption (\ref{h11}) into account, the proof of our theorem is completed once we prove the three following results:
	\begin{equation}
		\sup_{\theta\in \Theta_n}\sup_{\phi_{x,\theta}^{-1}\left(\frac{\rho_n k_{1,n}}{n}\right)\leq h\leq \phi_{x,\theta}^{-1}\left(\frac{ k_{2,n}}{\rho_n n}\right)}\left|\hat{r}_{h,\theta}(x)-r_{\theta_0}(x)\right|=O\left(f^{-1}\left(\frac{k_{2,n}}{\rho_n n}\right)^{\beta_0}\right)+O_{a.co.}\left(\sqrt{c_n}\right),\label{orden_kNN_theta_estimado}
	\end{equation}
	\begin{equation}
		\sum_{n=1}^{\infty}\sum_{m=1}^{n^{\alpha}}\sum_{k=k_{1,n}}^{k_{2,n}}\mathbb{P}\left(H_{k,x,\theta_m}\leq \phi_{x,\theta _m}^{-1}\left(\frac{\rho_n k_{1,n}}{n}\right)\right)< \infty, \label{inf_k1}
	\end{equation}
	and
	\begin{equation}
		\sum_{n=1}^{\infty}\sum_{m=1}^{n^{\alpha}}\sum_{k=k_{1,n}}^{k_{2,n}}\mathbb{P}\left(H_{k,x,\theta_m}\geq \phi_{x,\theta _m}^{-1}\left(\frac{k_{2,n}}{n\rho_n}\right)\right)< \infty,\label{sup_k2}
	\end{equation}
	where we have denoted $$c_n=\frac{\log n}{nf\left(\lambda f^{-1}(\rho_n k_{1,n}/n)\right)}.$$ On the one hand, the proof of (\ref{orden_kNN_theta_estimado}) is a direct consequence of Theorem \ref{teorema theta estimado}(a). Specifically, taking assumptions (\ref{h33bis}), (\ref{h33bisbis}) and (\ref{h888bis}) into account, it suffices to consider $a_n=\lambda f^{-1}\left(\frac{\rho_n k_{1,n}}{n}\right)$ and $b_n=\delta f^{-1}\left(\frac{ k_{2,n}}{\rho_n n}\right)$ in Theorem \ref{teorema theta estimado}(a). On the other hand, to prove (\ref{inf_k1}) and (\ref{sup_k2}) we will use Lemma \ref{lemaA.4}. By means of such lemma, and in a similar way as in \cite{kara}, one obtains:	
	\begin{equation}
		\mathbb{P}\left(H_{k,x,\theta}\leq \phi_{x,\theta _m}^{-1}\left(\frac{\rho_n k_{1,n}}{n}\right)\right) \leq \exp\left\{-\frac{(1-\rho_n) k_{1,n}}{4}\right\}   +   \exp\left\{-\frac{(1-\rho_n)^2 k_{1,n}}{4 \rho_n}\right\}      \label{exp1} 
	\end{equation}
	and
	\begin{equation}
		\mathbb{P}\left(H_{k,x,\theta}\geq \phi_{x,\theta _m}^{-1}\left(\frac{k_{2,n}}{n\rho_n}\right)\right)\leq \exp\left\{-\frac{(1-\rho_n)^2 k_{2,n}}{2\rho_n}\right\}.\label{exp2}
	\end{equation}
	Note that the first summand in the bound (\ref{exp1}) is related to the case where $\min\{\omega^2, \omega\}=\omega$ while the second one, which was forgotten in \cite{kara}, corresponds to the case where $\min\{\omega^2, \omega\}=\omega^2$ (we have denoted $\omega=k/(\rho_n k_{1,n}) -1$); for details on the role of $\omega$, see Lemma \ref{lemaA.4}. Therefore one has that
	\begin{eqnarray}
		\sum_{m=1}^{n^{\alpha}}\sum_{k=k_{1,n}}^{k_{2,n}} \mathbb{P}\left(H_{k,x,\theta_m}\leq \phi_{x,\theta _m}^{-1}\left(\frac{\rho_n k_{1,n}}{n}\right)\right)&\leq& n^{\alpha} k_{2,n} (n^{-\frac{1-\rho_n}{4}\frac{k_{1,n}}{\ln n}} + n^{-\frac{(1-\rho_n)^2}{4\rho_n}\frac{k_{1,n}}{\ln n}}) \nonumber \\ &\leq& n^{\alpha + 1 -\frac{1-\rho_n}{4}\frac{k_{1,n}}{\ln n}} + n^{\alpha + 1-\frac{(1-\rho_n)^2}{4\rho_n}\frac{k_{1,n}}{\ln n}}
		\label{fin1}
	\end{eqnarray}
	and
	\begin{eqnarray}
		\sum_{m=1}^{n^{\alpha}}\sum_{k=k_{1,n}}^{k_{2,n}} \mathbb{P}\left(H_{k,x,\theta_m}\geq \phi_{x,\theta _m}^{-1}\left(\frac{k_{2,n}}{n\rho_n}\right)\right) 
		\leq n^\alpha k_{2,n} n^{-\frac{(1-\rho_n)^2}{2\rho_n}\frac{k_{2,n}}{\ln n}} \leq  n^{\alpha + 1-\frac{(1-\rho_n)^2}{2\rho_n}\frac{k_{2,n}}{\ln n}}. \label{fin2}
	\end{eqnarray}
	Finally, from Assumption (\ref{hnew}) together with the bounds (\ref{fin1}) and (\ref{fin2}), we obtain (\ref{inf_k1}) and (\ref{sup_k2}). This completes the proof of the theorem.

	\subsection{Proof of Corollary \ref{cor}}
	It suffices to check that the assumptions used in Theorem \ref{teorema theta estimado}(b) hold and then to write the corresponding rate of convergence for the particular case considered in Corollary \ref{cor}.
	
	\subsection{Proof of Corollary \ref{final-cor}}
	Trivial. 
	
\end{document}